\def\draft{n}

\documentclass{amsart}
\usepackage{amssymb, amsmath, latexsym, amscd}
\usepackage{epsfig}

\def\printname#1{
	\if\draft y
		\smash{\makebox[0pt]{\hspace{-0.5in}
			\raisebox{8pt}{\tt\tiny #1}}}
	\fi }
\def\lbl#1{\label{#1}\printname{#1}}

\newtheorem{theorem}{Theorem}[section]
\newtheorem{corollary}[theorem]{Corollary}
\newtheorem{lemma}[theorem]{Lemma}

\newtheorem{proposition}[theorem]{Proposition}
\newtheorem{problem}[theorem]{Problem}
\newtheorem{conjecture}[theorem]{Conjecture}
\newtheorem{challenge}{Challenge}
\newtheorem{ktok}{Theorem~\ref{imply}}
\newtheorem{four}{Proposition~\ref{4D}}
\newtheorem{bound}{Proposition~\ref{bounding}}
\newtheorem{torsion}{Proposition~\ref{2-torsion}}

\theoremstyle{definition}
\newtheorem{definition}[theorem]{Definition}
\newtheorem{definitions}[theorem]{Definitions}
\newtheorem{remark}[theorem]{Remark}

\newcommand{\A}{\mathcal A}
\newcommand{\B}{\mathcal B}
\newcommand{\C}{\mathcal C}
\newcommand{\G}{\mathcal G}

\newcommand{\K}{\mathcal K}
\newcommand{\F}{\mathcal F}
\newcommand{\Z}{\mathbb Z}
\newcommand{\s}{\mathcal S}
\newcommand{\Q}{\mathbb Q}
\newcommand{\T}{\mathcal T}

\newcommand{\Prim}{\operatorname{Prim}}

\newcommand{\ra}{\longrightarrow}

\newcommand{\gc}{CC-set }
\newcommand{\dcs}{CC-sets }
\newcommand{\dc}{\gc}
\newcommand{\ccf}{CC-scheme }
\newcommand{\ccfs}{CC-schemes }
\renewcommand{\c}{\mathcal C}
\newcommand{\f}{\mathcal S}

\begin{document}

\title{Grope cobordism and Feynman diagrams}  
\author[J. Conant]{James Conant}
\address{Department of Mathematics\\
         Cornell University\\
         Ithaca, NY 14853-4201}
\email{jconant@math.cornell.edu}

\author[P. Teichner]{Peter Teichner}
\address{Department of Mathematics,
         University of California in San Diego,
         9500 Gilman Dr, La Jolla, CA 92093-0112}
\email{teichner@euclid.ucsd.edu}

\thanks{The first author was partially supported by NSF VIGRE grant DMS-9983660.
The second author was partially supported by NSF grant DMS-0072775 and the
Max-Planck Gesellschaft.}
\begin{abstract}
We explain how the usual
algebras of Feynman diagrams behave under the grope degree introduced in
\cite{ct}. We show that the Kontsevich integral rationally classifies grope
cobordisms of knots in 3-space when the ``class'' is used to organize gropes.
This implies that the grope cobordism equivalence relations are highly
nontrivial in dimension~$3$. We also show that the class is not a useful
organizing complexity in $4$ dimensions since only the Arf invariant survives.
In contrast, measuring gropes according to ``height'' does lead to very
interesting $4$-dimensional information
\cite{cot}. Finally, several low degree calculations are explained, in particular we
show that S-equivalence is the same relation as grope cobordism based on the
smallest tree with an internal vertex. 
\end{abstract}
\maketitle

\large

\section{Introduction} 
In \cite{ct} we introduced the notion of a {\em grope
cobordism} between two knots in 3-space, which places Vassiliev theory in a
natural topological context. Gropes
are certain 2-complexes built out of several surface stages, whose complexity
can be measured by either the {\em class} (corresponding to nilpotent groups)
or the {\em height} (corresponding to solvable groups). The analogy to group
theory arises by observing that a continuous map $\phi$ of a circle (into some
target space) represents a commutator in the fundamental group if and only if
it extends to a map of a surface (which is the simplest possible grope,
of class~2 and height~1). Similarly, $\phi$ represents an element in the $k$-th
term of the {\em lower central series} (respectively {\em derived series}) of
the fundamental group if and only if it extends to a continuous map of a grope
of class~$k$ (respectively height~$k$).

In knot theory, one replaces continuous maps of a circle by smooth 
embeddings of a circle into 3-space. Accordingly, one should study {\em
embeddings} of gropes into 3-space. More precisely, one obtains two sequences of
new geometric equivalence relations on the set of knot types by calling
two knots equivalent if they cobound an {\em embedded} grope (of a specified class
or height). 

It is the purpose of this paper to show that the invariants associated to grope
cobordism are extremely interesting.
Let $\K$ be the abelian monoid of knot types, i.\ e.\ isotopy classes of oriented knots
in 3-space (under connected sum). We proved in \cite{ct} that the quotients $\K/G_k:=$
$\K$ modulo grope cobordism of class $k$ in 3-space, are in fact a finitely generated
abelian groups. In Section~\ref{sec:GR} we start the investigation of these groups by
showing that there is an epimorphism
$$
\B^g_{<k} \twoheadrightarrow \K/G_k,
$$
where $\B^g_{<k}$ is the usual (primitive) diagram space known from the theory
of finite type invariants but {\em graded by the grope degree}. More
precisely, $\B^g_{<k}$ is the abelian group generated by connected
uni-trivalent graphs of grope degree~$i,\quad 1<i<k$, with at least one
univalent vertex and a cyclic ordering at each trivalent vertex. The
relations are the usual IHX and AS relations. The grope degree is the
Vassiliev degree (i.e.\ half the number of vertices) {\em plus} the first
Betti number of the graph. Observe that both relations preserve this new
degree. 

We then show in Section~\ref{sec:GR} that as in the usual theory of finite type
invariants, the above map has an inverse, the {\em Kontsevich integral}, after tensoring
with the rational numbers.  See Definition~\ref{def:Z} for the definition of
the Kontsevich integral $Z^g_{<k}$:

\begin{theorem}\lbl{B}
$Z^g_{<k}$ induces an isomorphism of $\Q$-vector spaces
$$
\K/G_k\otimes\Q\cong \B^g_{<k}\otimes \Q
$$
\end{theorem}
This result was inspired by the recent discovery of Garoufalidis and Rozansky \cite{gr}
that the Kontsevich integral not only preserves the Vassiliev filtration but also the
``loop filtration'', where one grades diagrams by the first Betti number (and
correspondingly the clasper surgeries are reorganized). We decided to give an
independent proof of Theorem~\ref{B}, by using properties of the Kontsevich integral
explained in \cite{Aa} as well as a result from \cite{ct} which says that grope
cobordisms can be refined into simple clasper surgeries.
In Conjecture~\ref{integral}, we take a guess at what the groups $\K/G_k$
could be integrally.  There is an analogue of this Theorem~\ref{B} which says 
that {\em capped} grope cobordism is rationally computed by 
$\B^v_{<k}$, which is the same diagram space as above but graded by the 
Vassiliev degree. This latter result follows by using work of Habiro, and was
announced in \cite{h2}. 

If one uses class as an organizational tool for {\em grope
concordance}, i.e.\ for gropes embedded in $S^3 \times [0,1]$ and with boundary in
$S^3 \times \{0,1\}$, then the theory collapses:

\begin{four} 
For each $k\geq 3$, two knots  $K_i\subset S^3 \times\{i\}$, are class $k$ grope
concordant if and only if their Arf invariants agree.
\end{four}

It should be mentioned that Schneiderman has independently given a beautiful
geometric argument for the above fact: He directly constructs a {\em weak
Whitney tower of class $k$} in 4-space, cobounding two knots with equal Arf
invariants. One can then turn this weak Whitney tower into a grope concordance
of class $k$.

It turns out that in order to
derive interesting information about knot concordance (i.e.\
$4$-dimensional knot theory), one needs to imitate the {\em derived series}
of a group geometrically. This can be done by restricting attention to
gropes which grow symmetrically from the root. Such gropes have a {\em
height} $h$ and the class $k$ can be calculated as
$k=2^h$, exactly as for group commutators. It was shown by Cochran, Orr and
Teichner in \cite{cot} that  {\em symmetric} grope concordance filters the
knot concordance group in such a way as to yield all known concordance
invariants in the first few steps, leaving a huge tower of concordance
invariants yet to be discovered. The first new graded quotient (above the
one leading to Casson-Gordon invariants) was shown to be nontrivial in
\cite{cot} by using certain von Neumann signatures associated to solvable quotients
of the knot group. It is now known \cite{CoT} that all the graded quotients are
nontrivial.

\begin{challenge}
Try to understand the equivalence relation of symmetric grope cobordism in
3-space. In particular, determine the role of the von Neumann signatures.
\end{challenge}

We will show that very interesting things happen even at very small
heights. At height $h=1.5$ one gets isomorphism of Blanchfield forms
(i.e.\ S-equivalence) for grope cobordisms in 3-space, whereas the
$4$-dimensional analogue gives {\em cobordism} of Blanchfield forms. So in
this setting the ``kernel'' from dimensions~$3$ to $4$ is given by connected
sums $K \# K^!$, where $K$ is any oriented knot and $K^!$ is its reversed
mirror image, the concordance inverse.

The reader might be irritated about the occurrence of the non-integral
height $h=1.5$ but that's just a special case of the following equivalence
relations on knots: Fix a rooted tree type $\T$ and consider only grope
cobordisms of type $\T$.  The notation $\T$ (respectively $c\T$) in the
following table refers to the equivalence relation given by grope
cobordisms (respectively capped grope cobordism) in
$3$-space using gropes of tree-type $\T$, as explained in \cite{ct}. One can
also study grope cobordism in $S^3 \times [0,1]$ which is denoted by $\T^4$
below. Note that in dimension~$4$ there is no difference between capped and
uncapped grope cobordism because intersections and self-intersections of the
caps can always be pushed down into the bottom stage. The following table
summarizes our calculations in Section~\ref{sec:low}.

\begin{theorem} \lbl{table}
For the smallest rooted tree types $\T$, the grope cobordism (respectively grope
concordance) relations are given by the following table:
\vspace{5mm}

\noindent\begin{tabular}{|r||c|c|c|}
\hline Tree Type $\T$ &
$\K/c\T$ &
$\K/\T$ &
$\K/\T^{4}$\\
\hline\hline
\epsfig{file=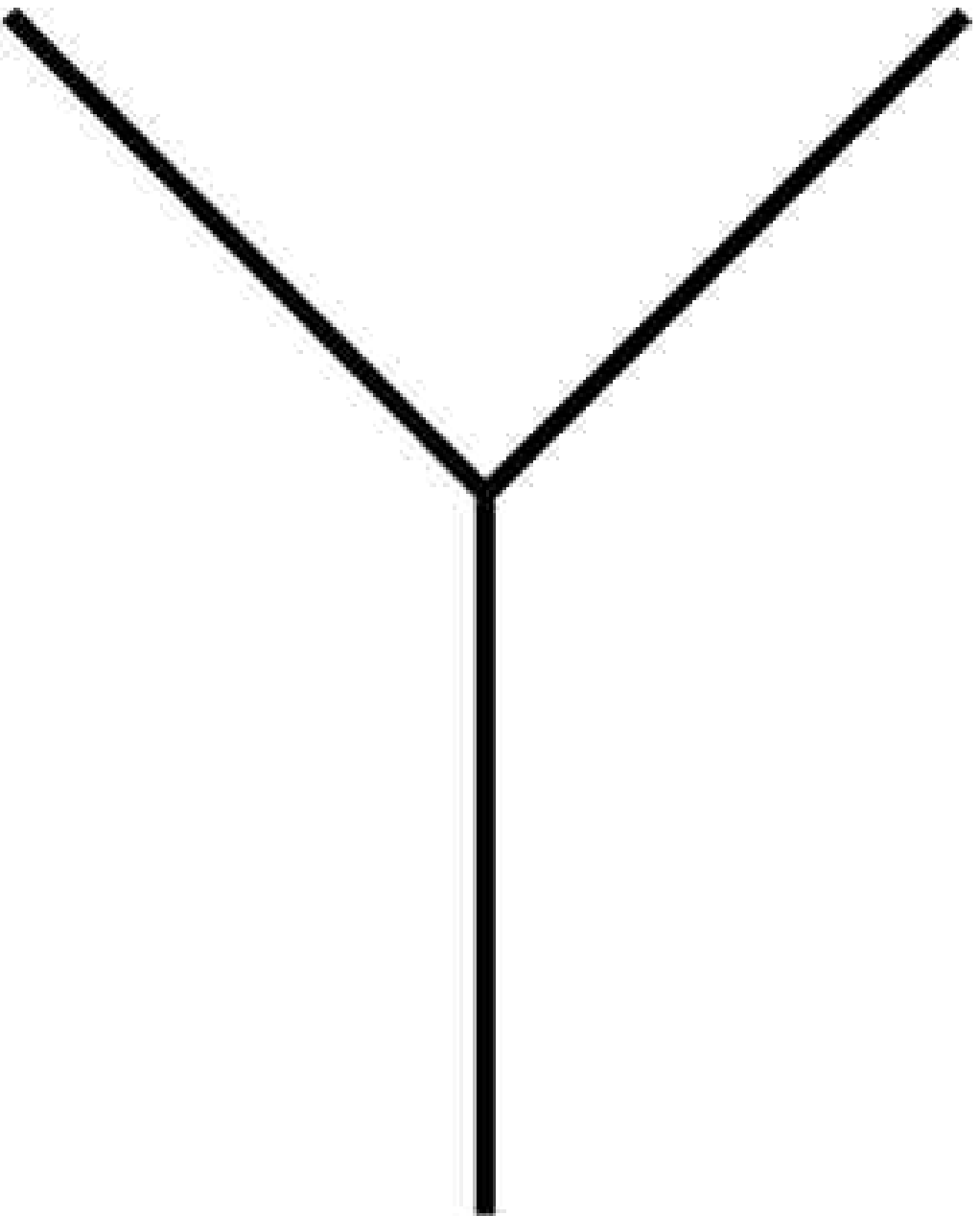,%
           height=.3cm} &
$\{0\}$&$\{0\}$&$\{0\}$\\
\hline
\epsfig{file=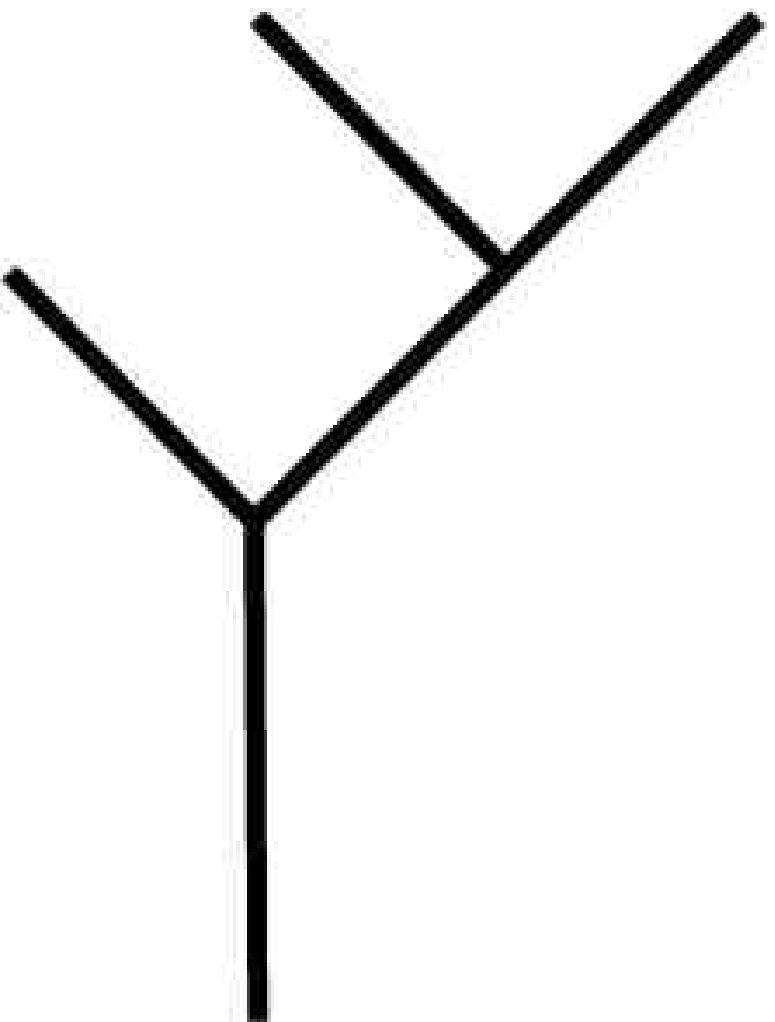,%
           height=.5cm} &
$\Z (c_2)$ &
$\Z/2 (\text{Arf})$ &
$\Z/2 (\text{Arf})$\\
\hline
\epsfig{file=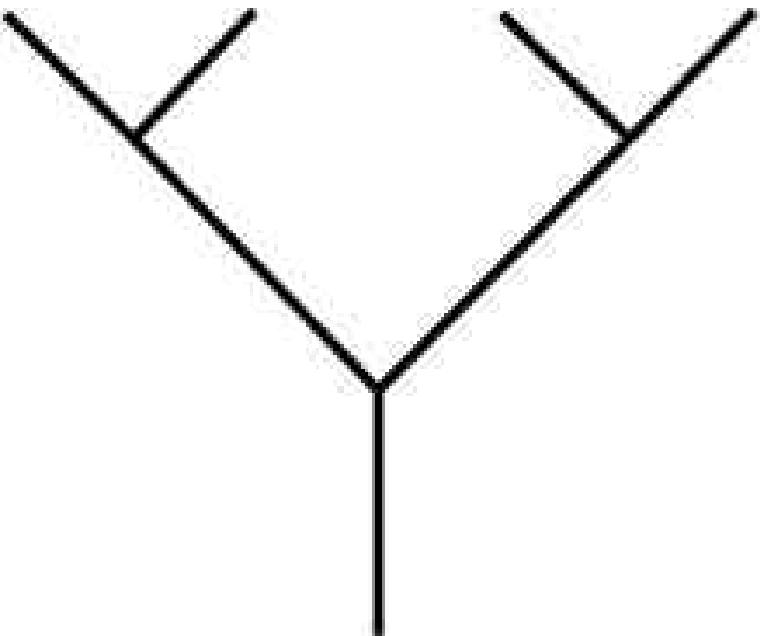,%
           height = .5cm} &
$\Z (c_3)\oplus \Z (c_2)$ &
$\Z (c_2)$ &
$\Z/2 (\text{Arf})$\\
\hline
\epsfig{file=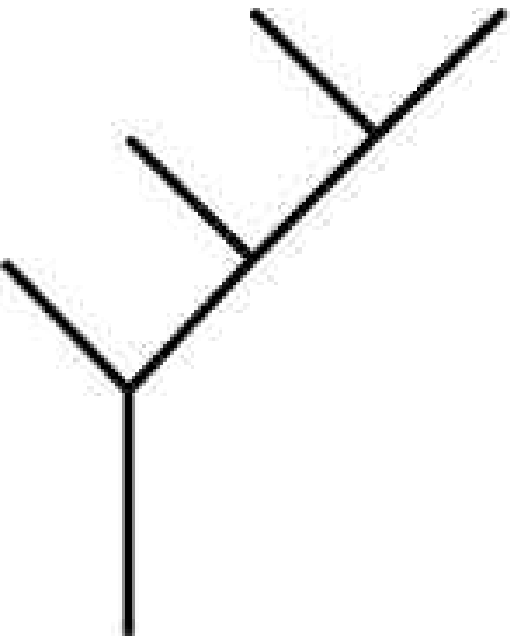,%
           height = .5cm} &
$\Z (c_3)\oplus \Z (c_2)$ &
$\Z (c_2)$ &
$\Z/2 (\text{Arf})$\\
\hline
\epsfig{file=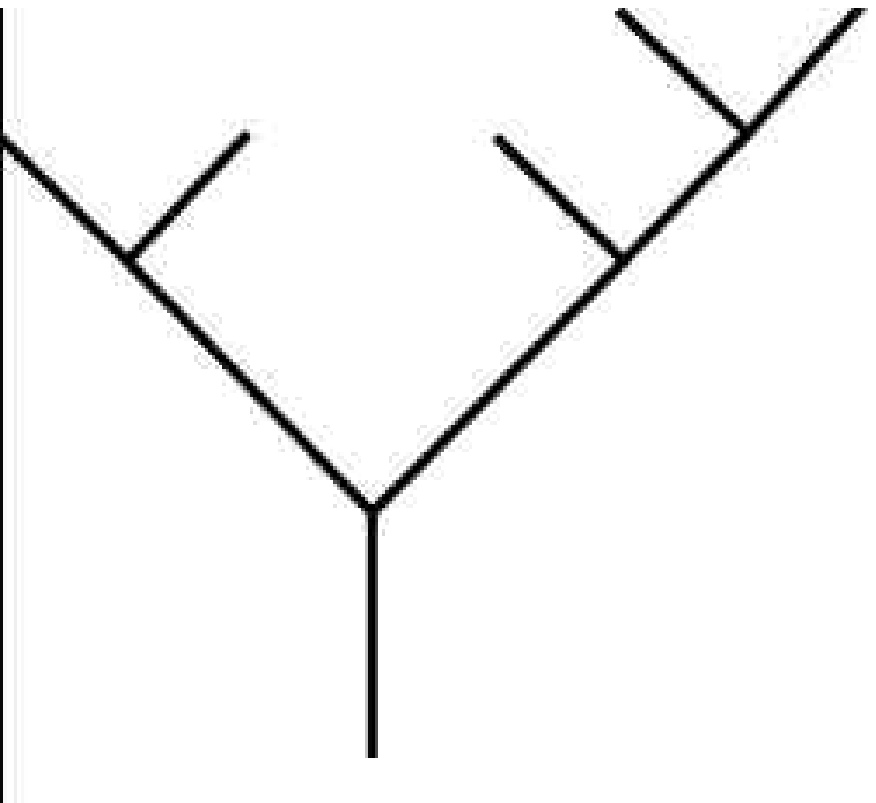,%
           height= .8cm} &
$\begin{array}{c}
\Z (c_4)\oplus \Z (c_4^\prime) \oplus\\
\Z (c_3)
\oplus \Z (c_2)
\end{array}$ &
$\Z/2 (c_3) \oplus \Z (c_2)$&
$\Z/2 (\text{Arf})$\\
\hline
\epsfig{file=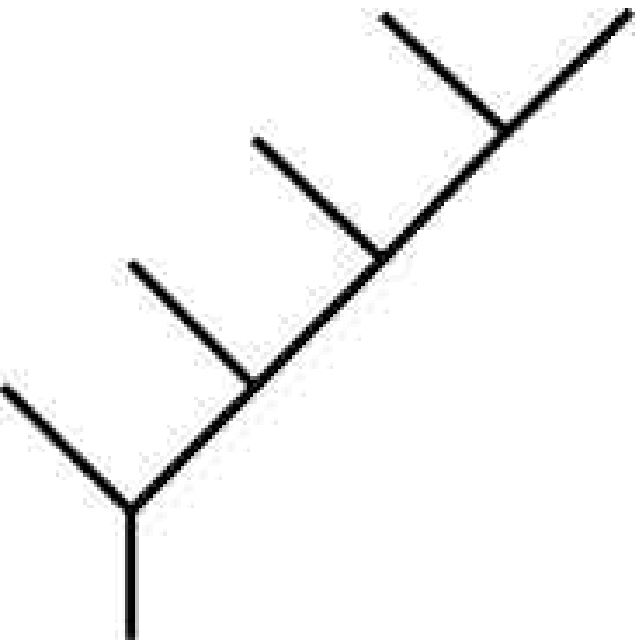,%
           height= .8cm} &
$\begin{array}{c}
\Z (c_4)\oplus \Z (c_4^\prime) \oplus\\
    \Z (c_3)
\oplus \Z (c_2)
\end{array}$ &
$\Z/2 (c_3) \oplus \Z (c_2)$&
$\Z/2 (\text{Arf})$\\
\hline
\epsfig{file=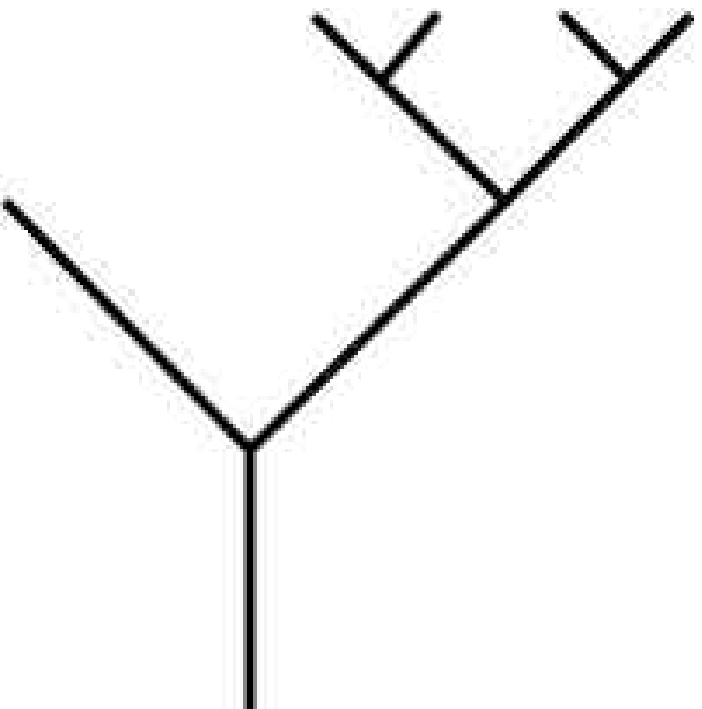,%
           height =.8cm} & {\bf ?}&
$\begin{array}{c}
\text{$S$-equivalence}\\
\text{or Bl-forms}
\end{array}
$ &
$\begin{array}{c}
\text{cobordism of}\\
\text{Bl-forms}
\end{array}
$\\
\hline
\end{tabular}
\vspace{5mm}

Here the $c_i$ are generating degree~$i$ Vassiliev invariants, and $Bl$ is
the Blanchfield form.
\end{theorem}

\begin{corollary}\lbl{cor:c3}
 $c_3$ modulo two is an S-equivalence invariant.
\end{corollary}

\begin{challenge}
Understand the monoids $\K/\T$ for more tree types $\T$.
\end{challenge}

Our paper ends with Section~\ref{sec:misc} where we have
collected several results that are relevant in the context of grope
cobordism.
Recall that if a knot $K$ cobounds a grope with the unknot $U$, then $K$ and $U$
might very well be linked in a nontrivial way. Thus it is a much stronger
condition on $K$ to assume that it is the boundary of a grope. For example, if $K$
bounds a grope of class~$3$ in $S^3$ then the Alexander polynomial vanishes. 
The following statement generalizes this vanishing result by using the
``null-filtration'' of \cite{gr}, explained in Section~\ref{sec:null}. 

\begin{bound} 
If a knot $K$ bounds an embedded
grope of class $k$ in a $3$-manifold $M$, then the pair $(M,K)$ is
$(k-3)$-null equivalent to $(M,U)$.
\end{bound}

The converse is not true for $k=3$: 
In $M=S^3$, knots which are null equivalent to the unknot ($k=3$ above), are
exactly knots with trivial Alexander polynomial. On the other hand, knots which
bound a grope of class~$3$ in $S^3$ have a Seifert surface such that the rank of
the Seifert form equals the genus of the surface. It is shown in
\cite{gt} that this {\em minimal Seifert rank} condition is much stronger than
having trivial Alexander polynomial.

We invesitgate the behavior of orientation reversal:
\begin{torsion} 
Let $\rho$ be the map reversing a knot's
orientation. Then for every knot $K$ in the $k$-th term $G_k$ of the grope
filtration of $\K$, one has 
$$
K\equiv (-1)^k\rho(K)\mod G_{k+1}.
$$
\end{torsion}

Our final result implies that the main theorem of \cite{ct} can now be phrased as
follows: Grope cobordisms of class $k$ in $S^3$ are in 1-1 correspondence with
sequences of simple clasper surgeries of grope degree {\em exactly} $k$. In
\cite{ct} we had to allow simple clasper surgeries of grope degree~$\geq k$.

\begin{ktok} 
A simple clasper surgery of grope degree~$(k+1)$ may be realized by a sequence
of simple clasper surgeries of grope degree~$k$.
\end{ktok}

\vspace{1mm}

\noindent {\em Acknowledgments}: It is a pleasure to thank Stavros
Garoufalidis and Kazuo Habiro for useful discussions.

\tableofcontents

\section{Gropes, claspers and diagrams}

\subsection{Basic notions} 

For the reader's convenience we recall some basic notions and results from
\cite{ct}. Given a knot $K$ and a disjoint rooted clasper
$C$ in a $3$-manifold $M$, one can construct a grope cobordism $G(K,C)$ between
$K$ and $K_C$ in $M$, where $K_C$ is the surgery of $K$ along
$C$. The ``root'' of the clasper is a {\em simple} leaf (i.e.\ a leaf with an
embedded cap intersecting $K$ once) which makes sure that the 
ambient $3$-manifold $M$ is
unchanged. {\em Caps} of a clasper $C$ are disjointly embedded disks with
interiors disjoint from $C$, which bound some of the leaves of $C$. The grope
cobordism $G(K,C)$ is {\em capped} if all the leaves of the clasper have
caps, which are allowed to intersect the knot. (If there is a single cap disjoint
from $K$ then $K_C$ is isotopic to $K$.)
Finally, in a {\em simple} clasper, all the leaves have caps and the knot
intersects each cap exactly once. Thus this notion only makes sense for the pair
$(K,C)$.

The class of $G(K,C)$ is given by the {\em grope degree} of the clasper. Recall
that a clasper
$C$ has an underlying uni-trivalent graph $\Gamma$ which is obtained by removing
the leaves and collapsing to the spine. The grope degree of $C$ is defined as the
Vassiliev degree plus the first Betti number of $\Gamma$. The Vassiliev degree is
one half the number of vertices of $\Gamma$.

The construction of $G(K,C)$ depends on a choice of $b_1(\Gamma)$ many ``cuts'', which
turn $\Gamma$ into a rooted tree, giving the precise grope type of $G(K,C)$. 
Note that each cut increases the Vassiliev degree by one, but leaves the grope
degree unchanged, as it should. Note also that a cut introduces a pair of
Hopf-linked leaves into the clasper, and hence the resulting grope cobordism
cannot have two disjoint caps at the corresponding tips. This explains why
 capped grope cobordism corresponds to the Vassiliev degree: cuts are not
allowed (since caps must be disjointly embedded), and for trees the two degrees
agree.

Let $\K$ denote the monoid of oriented knot types in 3-space (with respect to connected
sum), and let $\K/G_k$ be $\K$ modulo the equivalence relation of grope cobordism of
class $k$. It turns out that this is an abelian group. 
By Theorem 2 of \cite{ct}, a grope cobordism of class
$k$ corresponds to a sequence of simple clasper surgeries of grope degree~$\geq
k$. We will show in Theorem~\ref{imply} that  a grope cobordism of class $k$
also corresponds to a sequence of simple clasper surgeries of grope degree
{\em exactly equal to} $k$.

Let $\K/G^c_k$ denote the abelian group of oriented knot types
modulo capped grope cobordism of class $k$ in 3-space. Capped grope
cobordism of class $k$ coincides with Vassiliev degree~$k$ simple clasper
surgeries, and with Vassiliev
$(k-1)$-equivalence. By letting $G_k$ be the subset of knots which are
class~$k$ grope cobordant to the unknot (and similarly $G^c_k$ for capped
grope cobordisms) we can form the associated graded quotients
$$
\G_k := G_k / G_{k+1} \quad\text{ respectively }\quad \G^c_k := G^c_k / G^c_{k+1}
$$
These graded quotients will be related in Lemmas~\ref{Phi} and
\ref{Phig} to certain Feynman diagrams of a fixed degree.

Let $\T$ be a rooted tree type. Then $\K/\T$ is defined as
the monoid of knots modulo $\T$ grope cobordism. Define $\K/\T^c$ to be the
monoid of knots modulo capped $\T$ grope cobordism. A capped $\T$ grope
cobordism corresponds to a sequence of simple tree clasper surgeries of type
$\T$. As a consequence, $\K/\T^c$ depends only on the unrooted tree type.

\subsection{Feynman diagrams}

Let $\widetilde{\A}_k$ denote the free abelian group generated by connected
trivalent graphs with $2k$ vertices and one distinguished (oriented) cycle,
such that each trivalent vertex has a cyclic orientation.\footnote{ Since
trivalent vertices that lie on the distinguished cycle can be canonically
oriented, the convention in much of the literature is to not orient these
vertices.} The distinguished cycle is often called ``the outer circle,'' and
the rest is sometimes called ``the dashed part." 

Define $\A_k$ to be the quotient of $\widetilde{\A}_k$ by the usual IHX and
AS relations. The AS (antisymmetry) relation is a relation of the form
$G_1+G_2=0$, where the $G_i$ differ by a cyclic orientation at a given vertex.
 The IHX relation is pictured in Figure
\ref{1tihx}. If the distinguished cycle runs through the part of the
diagram involved in an IHX relation, the relation is called an STU
relation. As proven in \cite{bn}, STU relations generate all IHX relations.
\begin{figure}[ht] 
\begin{center}
\epsfig{file=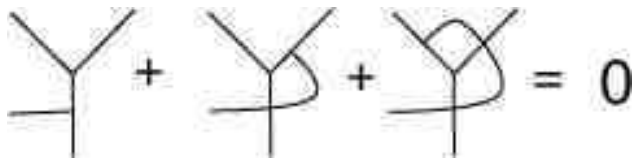,%
        height = 2cm}
\end{center}
\caption{The IHX relation.}\lbl{1tihx}
\end{figure}
Then $\A = \oplus_k \A_k$ is the well-known algebra of Feynman diagrams.
(Feynman diagrams enjoy a plethora of other names, but we will stick to this
one.)  The algebra structure is given by ``connected sum'', which turns out to
be well-defined. Moreover,
$\A$ is a graded Hopf algebra as explained in
\cite{bn}. The {\em primitive} elements $\Prim(\A)$ are generated by
diagrams which stay connected when removing the outer circle.
However, such diagrams are {\em not} closed under the STU relation, so it is
convenient to also consider the group
$\A^I :=\A_{>0} /\A_{>0}^2$ of {\em indecomposable} elements. Since $\A$ is
a commutative and cocommutative Hopf algebra, a famous theorem of Milnor and
Moore implies that, say over $\Q$, $\A$ is a polynomial algebra in the primitive
elements. In particular, 
$$
\Prim(\A)\otimes\Q \cong \A^I\otimes\Q.
$$
There is an analogous abelian group of diagrams $\B$ defined just like $\A$
but without an outer circle. Here the subspace of connected diagrams is
closed under IHX and AS relations, so that the above problem for primitives
does not occur. There is an averaging map $\chi:\B\to\A$ which puts back the
outer circle in all possible ways. It was shown in \cite{bn} that this is
rationally an isomorphism. 

Since we will only be interested in connected diagrams, we {\em define}
$\B^v_k$ to be the abelian group spanned by connected uni-trivalent graphs
with $2k$ vertices (such that each trivalent vertex has a cyclic orientation,
and with at least one univalent vertex) modulo the AS and IHX relations. Since
we are dealing with unframed knots, we just {\em set} the groups $\B^v_k:=0$ for
$k=0,1$. Define the graded abelian group
$$
\B^v:=\oplus_k\B^v_k.
$$
The superscript `$v$' indicates that we
are using the Vassiliev degree, and it also serves to distinguish from the
usual group $\B$ (which would also contain non-connected diagrams). As a
consequence of these definitions, one has an averaging isomorphism of graded
vector spaces
\begin{equation}\lbl{eq:average}
\chi:\B^v\otimes\Q \overset{\cong}{\ra} \A^I\otimes\Q.
\end{equation}
We shall show that these groups, modulo terms of degree~$>k$, are
isomorphic to $\K/G^c_k\otimes\Q$.

\subsection{Maps relating capped gropes and Feynman diagrams}

There is a well-known
map that sends a chord diagram to an alternating sum over crossing
changes on the unknot guided by the chords. 
This map can in fact be extended to all diagrams, as in the next lemma,
whose proof is found in the next section. 
 
\begin{lemma}\lbl{Phi}
For each $k>1$, there is an epimorphism
$$
\Phi_k: \A^I_k \longrightarrow \G^c_k.
$$
defined by sending a diagram to the alternating sum
over clasper surgeries corresponding to each connected component of ``the dashed
part."
\end{lemma}

The (unframed) Kontsevich integral is a map on isotopy classes of oriented
knots
$$
Z: \K \longrightarrow \widehat\A:= \prod_k \A_k\otimes\Q
$$ 
which sends connected sums of two knots to their {\em product} in
$\widehat\A$. The image of $Z$ lies inside the group-like elements of the
complete Hopf algebra $\widehat\A$. Thus we may compose it with the {\em
logarithm} in this complete Hopf algebra to end up in the subspace of
primitive elements. Then $\log Z$ takes connected sum to addition! 
Now decompose $\log Z$ according to Vassiliev degree and denote the degree~$k$
part by
$$
Z_k^v:\K \longrightarrow \Prim(\A_k)\otimes\Q \cong \A^I_k\otimes\Q.
$$
Then $Z_k^v$ factors through simple clasper surgeries of degree~$k+1$ by
\cite{h2}, and hence through $\K/G^c_{k+1}$ by \cite{ct}. Restricting to
knots in
$G^c_k$ we get a homomorphism
$$
Z_k^v: \G^c_k \longrightarrow \A^I_k\otimes\Q.
$$
The definition of $\Phi_k$ given in the next section, and in particular the
fact that it extends the original definition on chord diagrams, implies by
the universality of the Kontsevich integral the following

\begin{lemma}\lbl{Z}
$Z_k^v \circ \Phi_k = Id$.
\end{lemma}
This shows that $\Phi_k$ is an isomorphism modulo torsion because $\A_k$ are
finitely generated.

\subsection{Clasper forests and Feynman diagrams}
It is the goal of this section to show that the map $\Phi_k$ extends to
all diagrams in a natural way. This will lead to Lemma~\ref{difference}
which will be used in the proof of Theorem~\ref{calc}.
 This material is well-known to the experts and was essentially
announced by Habiro in \cite{h2}. A complication we face is that we
include all capped claspers, not just those which are trees. Since Habiro
usually restricts to trees, we need to include some extra arguments for
the general case.

We begin with some definitions which will be used in this section.
\begin{definitions}
\begin{itemize}
\item A \emph{\dc} is a union of finitely many disjoint capped claspers
on a knot $K$.
\item The degree of a \dc is the minimum of the degrees of each
connected component.
\item A \emph{\ccf} is a collection $\{\c_1,\c_2,\ldots, \c_l\}$ where
each
$\c_i$ is a \dc.
\item The degree of a \ccf  is the sum of the degrees of each
$\C_i$.
\item A \ccf is called \emph{simple} if each \dc is a single simple
clasper.
\item Let $\f =\{\c_1,\c_2,\ldots, \c_l\}$ be a \ccf on a knot $K$.
Define 
$[K;\f]
 = \sum_{\sigma\subset \f}(-1)^{l+|\sigma |}
K_\sigma
\in \Z [\K],
$ where $K_\sigma$ is the knot modified by each \dc in
$\sigma$.
\end{itemize}
\end{definitions}

We list a couple of useful facts about brackets.

\begin{lemma}\lbl{goussarovlem}
Let $\c_1, \c_2$ be disjoint \dcs on $K$ and denote by $\c_1\cup\c_2$ the
\dc which is the union of $\c_1$ and $\c_2$. If $\f$ is a
\ccf on $K$, disjoint from $\c_1\cup\c_2$, then
$$
[K; \{\c_1\cup\c_2\}\cup \f] = [K; \{\c_1\}\cup \f] + [K_{\c_1};
\{\c_2\}\cup
\f],
$$
\end{lemma}

See \cite{g},
Lemma~5.2 mutatis mutandis.

\begin{figure}[ht]
\begin{center}
\epsfig{file=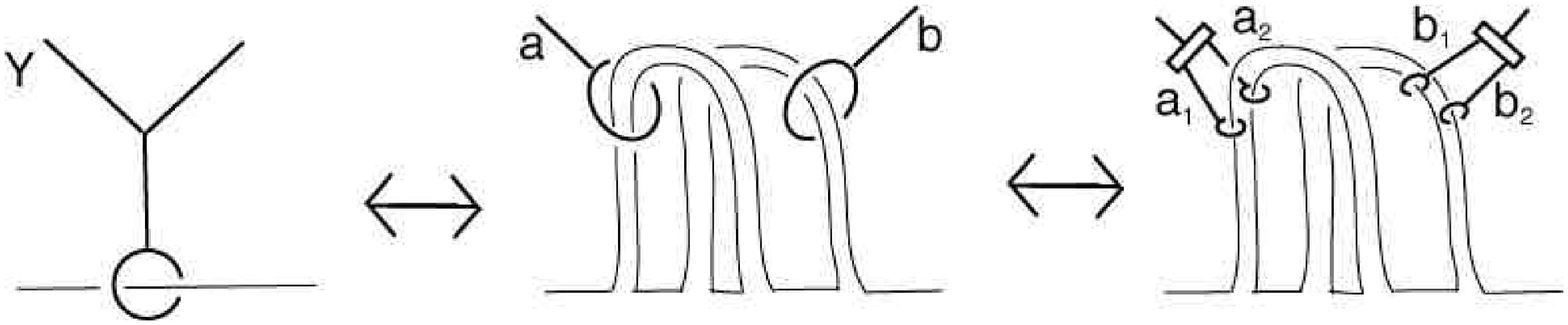,%
        height=3cm}
\end{center}
\caption{A clasper identity.}
\lbl{morse}
\end{figure}

\begin{lemma}\lbl{goodlem}
Consider the central part of Figure~\ref{morse}, 
with two disjoint capped claspers $a,b$, disjoint from a \ccf $\f$. Then
$$
[K;\{a\}\cup \{b\}\cup \f] = [K ; \{a \cup b\}\cup \f]
$$
\end{lemma}
\begin{proof}
If $\sigma\subset \f$ then
$K_{\{a\} \cup\sigma} = K_\sigma=K_{\{b\}\cup\sigma}$. 
This implies $[K_a;\f]=[K_b;\f]=[K;\f]$ and hence
\begin{align*} 
[K;\{a\}\cup \{b\}\cup \f] &= [K;\f] -[K_b;\f] - [K_a;\f] + [K_{a\cup
b};\f]\\  &= -[K;\f]+[K_{a\cup b};\f]\\
&=[K; \{a\cup b\}\cup \f]
\end{align*}
\end{proof}

The linear span of all knot types, $\Z[\K]$, can be filtered by
defining
$\F^v_k\subset \Z[\K]$ to be the linear span of all brackets $[K;\f]$
where $\f$ is of degree~$k$. 

\begin{lemma}\lbl{cong} 
If $C$ is a clasper on a knot disjoint from a
\ccf $\f$, then $[K;\f] - [K_C;\f] \in \F^v_{|\f|+|C|}$.
\end{lemma}

\begin{proof}
This is just the simple calculation $[K;\f\cup \{C\} ] = [K; \f] - [K_C;
\f]$. 
\end{proof}

The following proposition explains the superscript ``$v$."
\begin{proposition}\lbl{newprop}
$\F^v_k\subset \Z[\K]$ agrees with the usual Vassiliev filtration of
$\Z[\K]$.
\end{proposition}
\begin{proof}
One can take the usual Vassiliev filtration to be defined by brackets
where each \dc is a single simple degree $1$ clasper. 
Thus $\F^v_k$ contains the usual Vassiliev filtration.

The converse follow from Theorem 6.7 (3) of \cite{h2}, with one complication.
Namely, that Theorem is phrased only for tree claspers, but it can be easily
enhanced to work for all capped claspers.
The proof works by applying the move from Lemma~\ref{goodlem}, thereby breaking
a clasper into two claspers of lower degree. Iterating, one reduces all of the trees to
degree one capped claspers. This fails for graphs if in Figure~\ref{morse}, the two claspers
$a$ and $b$ are the same clasper. In this case we've still made progress since there
are a reduced number of loops. To be precise, induct on the number of internal vertices.
(Compare the proof of Lemma~\ref{separate}.)
\end{proof}

\begin{lemma}\lbl{separate}
On a knot $K$, let $\f_1=\{ C^1_1,\ldots, C_l^1\}$ and
$\f_2=\{C^2_1,\ldots,C^2_l\}$ be two \ccfs of degree
$k$, where each \dc consists of a single clasper.
Suppose the \ccfs differ only by a homotopy of an edge of one of the
claspers.
 Then
$$
[K;\f_1] - [K;\f_2]\in\F^v_{k+1}.
$$
\end{lemma}
\begin{proof}
 A homotopy of an edge can be
realized by a degree one clasper, one leaf of which links the edge and the other leaf being
embedded arbitrarily away from the caps. Using the zip construction, this can be realized by
degree one claspers $E$ where one leaf links the edge, and the other leaf either: is a meridian
to an edge, is a meridian to the knot, or is a trivial 1-framed leaf.
This last case is covered by Corollary~\ref{frame}. We prove the other two cases by induction
on the number of trivalent vertices.
The base case is when all the claspers $C^i_j$ are of degree one (eyeglasses). Then the leaves
of the clasper
$E$ can be slid off the end of the (up to two) eyeglasses that it links, introducing two new
intersections with the knot for each one. Let the slid clasper be called $E^\prime$. Then
we have argued
$$[K;\f_2]=[K_{E^\prime};\f_1].$$ By Lemma~\ref{cong} we are done. 
Now suppose that a pair $C^1_j,C^2_j$ has a trivalent vertex. 
Apply Figure~\ref{morse} to a leg of $C^i_j$, denoting the resulting
claspers by $a^i,b^i$. (If $E$ happens to link the leg we are expanding, push it off so that
it links the knot instead.)
  If $a^i$ and $b^i$ are different, then Lemma~\ref{goodlem}
implies that, defining $\f'_i= \{C^i_1,C^i_2,\ldots,C^i_{j-1},
a^i,b^i,C^i_{j+1},\ldots,C^i_l\}$,
$$[K;\f_i]=[K;\f'_i]$$
and since $\f'_1,\f'_2$ differ by the homotopy coming from $E$, 
$[K;\f'_1]=[K;\f'_2] \mod \F^v_{k+1}$ by induction.
If $a^i=b^i$, life is even simpler: define $$\f''= \{C^i_1,C^i_2,\ldots,C^i_{j-1},
a^i,C^i_{j+1},\ldots,C^i_l\}.$$ Then $[K;\f]=[K;\f'']$ and we are done by induction as above.
\end{proof}

\begin{lemma}\lbl{preinverse}
Let $C$ be a tree clasper on a knot $K$. There is a clasper $\widetilde{C}$ in a regular
neighborhood of $C$ such that $K_{C\cup\widetilde{C}}=K$. Moreover, the leaves of $C$ and
$\widetilde{C}$ are parallel (in the given framing of $C$) and one may arrange the edges
of $\widetilde C$ which go through a regular neighborhood of a leaf of $C$ to be
parallel to the leaves of $C$.
\end{lemma}
In Figure~\ref{eyeglasses}(a), the above process is depicted in the vicinity
of a Hopf-pair of $C$. (The two degree $1$ claspers labelled
$E$ are not relevant now.) 
\begin{proof}
The clasper $\widetilde{C}$ is constructed using Figure~27 of \cite{h2} and the version
of the zip construction in \cite[Sec.4]{ct}.
\end{proof}

\begin{lemma}\lbl{inverse}
Let $C$ and $C^\prime$ be two simple claspers of degree $k$ on a knot $K$,
which differ by a single half-twist. Then $K_{C} + K_{C^\prime} - 2K \in
\F^v_{k+1}$.
\end{lemma}
\begin{proof}
Insert Hopf-pairs of leaves to make $C$ a tree clasper and let $\widetilde{C}$ be an inverse
to $C$, as in Lemma~\ref{preinverse}. For each (non-root) cap  of $C$ and
$\widetilde{C}$ there is a degree~1 clasper (or eyeglass), surgery on which pushes
the knot out of the cap. Choose one of the leaves in each Hopf pair of
$C$. As this can be regarded as an additional cap of the tree clasper, there is an
eyeglass which pushes everything out of this leaf, just like above. See
Figure~\ref{eyeglasses} for an example of these eyeglasses $E_i$.  
Let $\mathcal E_1$
denote the \ccf whose elements are the (at least
$k$) eyeglasses which push things out of $C$, and let $\mathcal E_2$ be the
same, only for $\widetilde{C}$. The claspers in $\mathcal E=\mathcal
E_1\cup\mathcal E_2$ are capped, when considered on the knot
$K=K_{C\cup\widetilde{C}}$. Note that $K_\sigma = K$ if $\sigma$ has
nontrivial intersection with $\mathcal E_1$ and with $\mathcal E_2$. 
If $\emptyset\neq\sigma\subset\mathcal E_1$, then
$K_\sigma=K_{\widetilde{C}}$. If $\emptyset\neq\sigma\subset\mathcal E_2$
then $K_\sigma=K_{C}$. This analysis implies that
$$-[K;\mathcal E] = K_{C} + K_{\widetilde{C}}-2K.$$ 
The left-hand side is in $\F^v_{2k}\subset \F^v_{k+1}$. 

\begin{figure}[ht]
\begin{center}
\epsfig{file=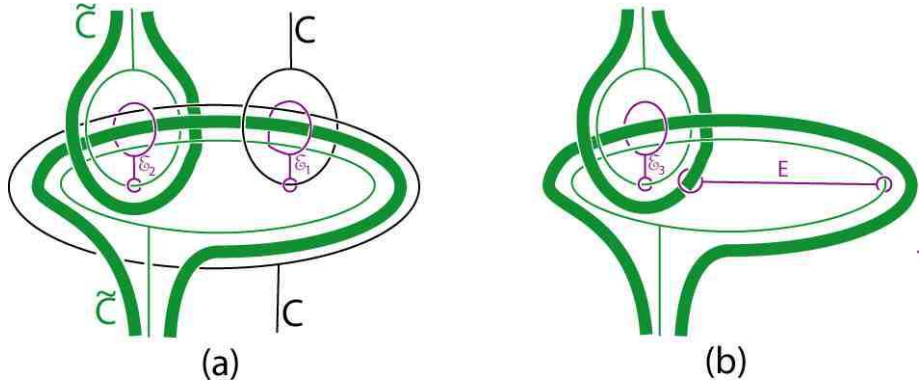, height=5cm}
\caption{From the proof of Lemma~\ref{inverse}. The thick
lines are where edges of $\widetilde{C}$ might be.}\lbl{eyeglasses}
\end{center}
\end{figure}

This is almost what we want, except that $\widetilde{C}$ has edges that wind
around itself and may wind through the Hopf pairs of leaves. 
Let $\mathcal E_3$ be the \ccf consisting of $k$ degree $1$ claspers on
$K_{\widetilde{C}}$, which push things out of the leaves.
There is an additional clasper, denoted by $E$ which pushes the edges
out of the other disk in each Hopf pair. It is a union of eyeglasses, one for each
leaf, see Figure~\ref{eyeglasses}(b). Now 
$$
[K_{\widetilde{C}}; \mathcal E\cup\{ E\}]=
K_{\widetilde{C}}-K_{\overline{C}},
$$ 
where $\overline{C}$ is formed from
$\widetilde{C}$ by applying $E$.
Symmetrically we can push the edges out of the other leaves in each Hopf
pair, modulo $\F^v_{k+1}$ to get the clasper
$C^\prime$ in the statement of the lemma.
\end{proof}

\begin{corollary}\lbl{frame}
Let $C_1$ and $C_2$ be two claspers on a knot $K$ which differ by a full
twist along an edge. Then $K_{C_1}-K_{C_2}\in\F^v_{k+1}$.
\end{corollary}
\begin{proof}
Let $C$ be a clasper that differs from both $C_1$ and from $C_2$ by a
half-twist. Then Lemma~\ref{inverse} implies that $K_{C_1}=2K-K_{C}
=K_{C_2} \mod \F^v_{k+1}$. 
\end{proof}

Let $\mathcal{CF}_k$ be the set of pairs $(K,\f)$ where $K\in\K$ and $\f$
is a simple \ccf of degree~$k$ on the knot $K$. These are
considered up to isotopy.  The bracket defines a
map $\mathcal{CF}_k \to \Z[\K ]$. 
One can also define
a map
$$
\Psi: \mathcal{CF}_k\to \A_k
$$
as follows: The unoriented trivalent graph, $D$,
associated to $(K,\f)$ is gotten by collapsing simple leaves to points and
forgetting the embedding. The distinguished cycle is the knot. The
orientation of the Feynman diagram is more subtle. 
Recall that each clasper $C_i$ is a certain thickening of a graph.
Hence $C_i$ may not be orientable (even though the thickened leaves are because they
have integral framing, and hence a full number of twists). However, each
$C_i$ may be made orientable by replacing $l_i$ bands by half-twisted bands, where
$l_i$ is the dimension of the $2$-torsion subgroup of
$H_1(C_i)$. Do this, and let $l=\sum l_i$ be the total number of half
twisted bands glued in to  every $C_i$. Now choose orientations of each new
thickening $C_i$. These give rise to cyclic orders of the in-coming edges at each
node. Denote the induced orientation of
$D$ by $or$. The choice of orientation also determines a normal direction.
Let
$m$ be the number of leaves where the knot pierces contrary to the
distinguished normal direction. Now $\Psi (K,\f) = (-1)^{m+l}(D,or)$.
$\Psi$ is well defined. It does not depend on the choice of orientations,
because the parity of the number of leaves of each $C_i$ is the same as the
parity of the number of nodes. The key thing to notice about this
orientation convention is that it switches under the introduction of a
half-twist in any edge of a clasper in $\f$.

Now we can define the homomorphism
$$
\Phi_k : \A^I_k \to \G^c_k 
$$
from Lemma~\ref{Phi} as follows:
For an oriented diagram $D\in \widetilde{A}_k$, define $\Phi_k(D):= cs_k([U;\f])$, for
some framed embedding $\f$ of the diagram $D$ on the unknot $U$, so that
$\Psi (U,\f) = D$. (Note that framed embedded diagrams can be regarded as simple \ccfs). 
Here
$cs_k:\Z[\K]\to \K/G^c_{k+1}$ is the map given by sending addition to connected sum.

By definition, $\Phi_k$ vanishes on diagrams of degree~$>k$ and, by Proposition~\ref{newprop}
lands in
$\G^c_k=\F^v_k/\F^v_{k+1}$. This, together with Lemma~\ref{separate} and
Corollary~\ref{frame}, implies that, extending linearly, we get a well defined map:
\begin{gather*}
\Phi_k\colon\widetilde{A}_k\to \G^c_k 
\end{gather*} 
This means that the choice of the framed embedding of the diagram $D$ is irrelevant.
Clearly $\Phi_k=0$ on separated diagrams, so it will factor through $\A^I_k$. 
Thus to get a map on $\A_k$, we only need to show that $\Phi_k$ vanishes on the STU- and
AS-relations.

\begin{lemma}
$\Phi_k$ vanishes on all STU and AS relations.
\end{lemma}
\begin{proof}
First we show it vanishes on AS relations. Suppose $D_1$ and $D_2$ are two
diagrams that differ by a cyclic order at one vertex. Let $C_1$ be a
clasper representing $D_1$: $\Phi_k(D_1)=U_{C_1}$. Let $C_2$ be the same,
except for three half-twists on the edges incident to the vertex. By our
orientation conventions, $\Phi_k(D_2)=U_{C_2}$. Applying
Lemma~\ref{inverse} three times, with $K=U$, we see that
$$
\Phi_k(D_1+D_2)=U_{C_1}\#U_{C_2}=0\quad\in\G^c_k,
$$
 as desired.  
Next we consider the STU relation.
Let $\widetilde{\Phi}_k :\widetilde{A}_k\to\Z[\K ]$ be the lift of $\Phi_k$ discussed
above, i.\ e.\ before applying the summation maps $cs_k$.
Consider an STU relation $D_s=D_t-D_u$, where $D_s$ has one more trivalent vertex than $D_t$
and
$D_u$.
Let Y be the component of the dashed part of $D_s$ which has the additional vertex. Then
\begin{gather*}
\widetilde{\Phi}_k (D_s) = [U; \{Y, C_2,\ldots ,C_l\}]
\end{gather*} where the clasper representing Y is also called Y, and is
pictured  in Figure \ref{morse}. The claspers $a$ and $b$ from Figure
\ref{morse}
 can be further subdivided into claspers $a_1,a_2,b_1,b_2$ using the zip
construction, assuming that
$a$ and $b$ are not the same clasper. By the previous lemmas we have

\begin{eqnarray*} & & [U;\{a\cup b, C_2,\ldots, C_l\}]\\ &=& [U;
\{a,b,C_2,\ldots,C_l\}]\\ &=& [U; \{a_1\cup a_2, b_1\cup b_2, C_2,\ldots,
C_l\}]\\ &=& [U;\{ a_1, b_1, C_2,\ldots, C_l\}] + [U_{b_1};\{a_1,b_2,C_2,\ldots
,C_l\}]\\ & & +[U_{a_1};\{ a_2,b_1,C_2,\ldots, C_l\}] + [U_{a_1\cup b_1};
\{a_2,b_2, C_2,\ldots, C_l\}]\\ &\equiv& [U;\{ a_1, b_1, C_2,\ldots, C_l\}] +
[U;\{a_1,b_2,C_2,\ldots ,C_l\}]\\ & & +[U;\{a_2,b_1,C_2,\ldots, C_l\}] + [U;
\{a_2,b_2, C_2,\ldots, C_l\}]
\mod \F^v_{k+1}(\Z\K)\\ &=& -\widetilde{\Phi}_k(D_t) +\widetilde{\Phi}_k(D_t)
-\widetilde{\Phi}_k(D_u) -
\widetilde{\Phi}_k(D_t) = \widetilde{\Phi}_k(D_t) - \widetilde{\Phi}_k(D_u)
\end{eqnarray*}

If, in Figure \ref{morse}, the two claspers $a$ and $b$ are really two ends
of the same clasper, we use the construction of Proposition 4.6 of
\cite{h2} instead. See Figure \ref{prop4.6} which is a clasper identity of
\cite{h2}. 
\begin{figure}[ht] 
\begin{center}
\epsfig{file=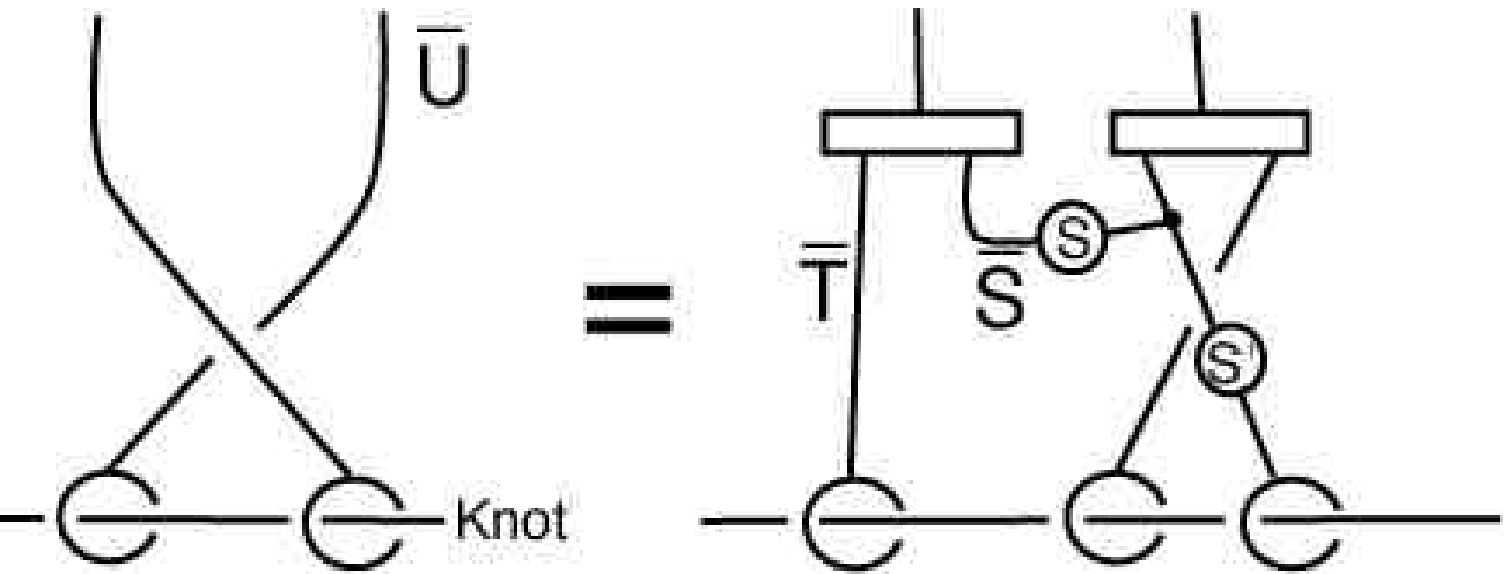,%
        height = 4cm}
\end{center}
\caption{A clasper identity.}\lbl{prop4.6}
\end{figure} 

Let $D_s,D_t,D_u$ be
the three diagrams in an STU relation such that the corresponding
$a$ and $b$ claspers are part of the same clasper. Then 
\begin{eqnarray*}  \widetilde{\Phi}_k (D_u) &=& [U;\{ \overline{U},
C_2,\ldots, C_l\}] = [U;
\{\{\overline{T},\overline{S}\}, C_2,\ldots, C_l\}]\\  &=& [U;\{ \overline{T},
C_2,\ldots, C_l\}] + [U_{\overline{T}};\{\overline{S}, C_2,\ldots, C_l\}]\\
&\equiv& [U;\{\overline{T},C_2,\ldots, C_l\}]+ [U;\{\overline{S}, C_2,\ldots,
C_l\}]\\ &=&
\widetilde{\Phi}_k (D_t) +
\widetilde{\Phi}_k(D_s)
\end{eqnarray*}
\end{proof}
The proof of Lemma~\ref{Phi} is now complete.

\begin{lemma}\lbl{difference}
Let $C$ be a simple clasper surgery of Vassiliev degree~$k$ on a knot $K$. Let $D$
be the diagram obtained by thinking of the clasper as the dashed part, and of the
knot as the outer circle. Then $\Phi_k(\pm D) = K_C\# K^{-1}\in \G^c_k.$
\end{lemma}

\begin{remark}
Here is a sketch of a proof that Lemma~\ref{difference} is true
rationally. The proof that it is true integrally is given below. One can
think of the Kontsevich integral of $K_C-K$ as the Aarhus integral of the
difference of two links
$K\cup L_C-K$, where $L_C$ is the link associated to the simple clasper. Now, it
is not hard to show that the lowest degree term of the Aarhus integral is exactly
the diagram corresponding to the graph type of the clasper. (By arguments
analogous to those in section 3.3 of \cite{gr}. In the absence of a knot $K$,
this is the statement that the Aarhus/LMO invariant is universal with respect to
Goussarov's Y-filtration.) Now the result follows by Lemma~\ref{Z}. 
\end{remark}

\begin{proof}[Proof of \ref{difference}.]
Let $\widetilde{C}$ be a clasper on the unknot with diagram $D$. Put in
half-twists so that it has the same sign as $C$.
 By definition,
$\Phi_k(D)=U_{\widetilde{C}}=U_{\widetilde{C}}\# K\# K^{-1}$. By
Lemma~\ref{separate}, we can move $\widetilde{C}$ into a position
corresponding to the clasper $C$ on the knot $K$. Using
Corollary~\ref{frame}, we can add an even number of half-twists until the
framings agree.
\end{proof}

\section{The grope degree}
We first introduce a second degree on the graphs $\Gamma$ generating $\B^v$:

\begin{definition}
The {\em grope degree} of $\Gamma$ is
$$
g(\Gamma):=b_1(\Gamma )+v(\Gamma). 
$$
Let $\B^g_k$ be the grope degree~$k$ part of the group $\B^v$.
\end{definition}
We note that the grope degree is preserved by the IHX and AS relations, and
hence really gives a new {\em graded} abelian group 
$$
\B^g:=\oplus_k \B^g_k.
$$

\subsection{Feynman diagrams and the grope degree}\lbl{sec:GR}
In a similar fashion to $\Phi_k$, we define a surjective map 
$$
\Phi^g_k : \B^g_k\to\G_k
$$
on connected Feynman diagrams of grope degree~$k$ as follows. For a connected
diagram $D\in\widetilde{\B}_k$, let $\widetilde{D}\in\widetilde{A}_k$ be a diagram formed by
attaching the univalent vertices of
$D$ to the outer circle in some order. Let $(U,C)\in\C\F^v_k$ satisfy
$\Psi (U,C) = \widetilde{D}$, where $C$ is a single clasper of grope degree~$k$, and $U$ is the unknot. Now define 
$$
\Phi_k^g(D) := U_C\in\K/G_{k+1}.
$$
 In the previous section we used various lemmas of Habiro \cite{h2} (for moving
claspers around modulo higher Vassiliev degree) in order to show that $\Phi_k$
is well defined. Section \ref{clgrdeg} contains the relevant lemmas in the case of
grope degree.
In fact, we need to show that $\Phi_k^g$ does not depend on the order in which we
attached the univalent vertices, and that AS and IHX get killed. The independence
of order follows from Lemma~\ref{lemmastar}(a), and the AS and IHX relations
follow from
\ref{lemmastar2}(a) and \ref{lemmastar2}(c).  
We have thus proven the following analogue of Lemma~\ref{Phi}.

\begin{lemma}\lbl{Phig}
For each $k>1$, there is an epimorphism
$$
\Phi^g_k: \B^g_k \longrightarrow \G_k
$$
\end{lemma}
In order to show that $\Phi^g_k$ is rationally an isomorphism, we want to
study the behavior of the Kontsevich integral with respect to the grope
degree. 

Recall that there is an isomorphism $\hat{\A}\cong\hat{\B}$ of graded {\em
algebras} given by the composition $\partial_\Omega \sigma$, where
$\sigma:\hat{\A}\to\hat{\B}$ is the inverse of the averaging map $\chi$ from
equation~\ref{eq:average}, and
$\partial_\Omega$ is the ``wheeling'' automorphism of $\hat\B$ \cite{thurston}. 
Here $\hat\A$ has the multiplication given by connected sum, whereas on $\hat\B$
the multiplication is just disjoint union. It follows that 
\begin{equation}\lbl{eq:log}
\log_{\hat\B}(\partial_\Omega  \sigma)= (\partial_\Omega 
\sigma)\log_{\hat\A}.
\end{equation}
Following \cite{gr}, define the {\em Euler degree} of a diagram in $\B$ to be the
number of {\em internal} trivalent vertices, which by definition are trivalent
vertices not adjacent to univalent vertices. 
It is called Euler degree because for uni-trivalent graphs $\Gamma$ one has
$$
e(\Gamma)=2(b_1(\Gamma)-b_0(\Gamma)).
$$
 Decompose the composition  $\log_{\hat\B}( \partial_\Omega  \sigma) Z$ according
to Euler degree to obtain
$$
Z^e_k:\K\to \B^e_k
$$ 
as the Euler degree~$k$ part of the ``Kontsevich integral''.
The following proposition follows by work of Garoufalidis and Rozansky
\cite{gr} using \cite{Aa}, but for the sake of completeness we provide an
argument using only \cite{Aa}.

\begin{proposition}
Let $C$ be a simple clasper of Euler degree~$n$ on a knot $K$. Then
$Z^e_m(K)-Z_m^e(K_C)=0$ for all $m<n$.
\end{proposition}

\begin{proof}
We factor through the Aarhus integral, $A$, of pairs $(M,K)$, \cite{Aa}. When
 normalized appropriately, $A(S^3,K)=\sigma Z(K)$. (In the definition of the
Aarhus integral we apply $\sigma$ to the knot $K$ as well as to the surgery
link components.) 
 We will assume familiarity with the
Aarhus integral in the following proof. 

Break $C$ into a union of $Y$'s.
If a $Y$ has no leaf that links the knot, it will be called \emph{internal} since
the corresponding trivalent vertex is internal.
 Each $Y$ has an associated $6$ component link. 
Let the three components in the Borromean rings be called the ``B" components,
and the three other components be called the ``L'' components. 
Let the link corresponding to the union of all Y's be called $L_C$. Fix balls
which meet the link at the Borromean rings corresponding to \emph{internal}
Y's as on the left hand side of Figure~\ref{borromean}. 
\begin{figure}[ht] 
\begin{center}
\epsfig{file=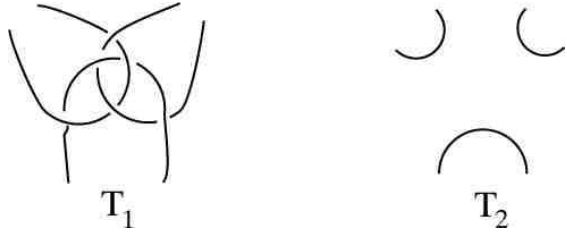,%
          height = 3cm}
\end{center}
\caption{The tangles $T_1$ and $T_2$. The difference $\mu:=Z(T_1)-Z(T_2)$ is
comprised of terms with at least one trivalent
vertex.}\lbl{borromean} 
\end{figure}
For each ball, there is an associated move that
replaces the tangle on the left of Figure~\ref{borromean} with the trivial one on
the right. Let
$S$ be the set of such moves, one for each internal vertex of $C$. Consider the
alternating sum $[K\cup L_C; S]$. Then the (alternating sum of) Kontsevich integrals
$Z[K\cup L_C;S]$  can be computed as the Kontsevich integral of the
difference of tangles in each ball, called $\mu$, glued to the
Kontsevich integral of the exterior to the balls. 
 See Figure~\ref{amalgamate}.
\begin{figure}[ht] 
\begin{center}
\epsfig{file=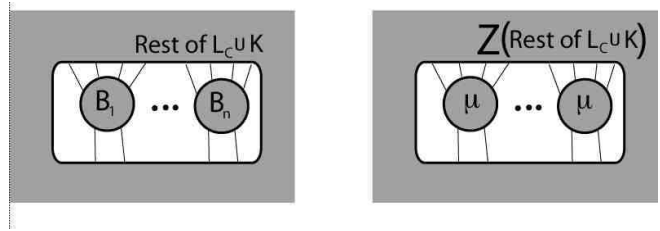,%
         height=3cm}
\end{center}
\caption{The Kontsevich integral of $[K\cup L_C;S]$ is computed by
gluing the Kontsevich integral at each ball $B_i$ to the
Kontsevich integral of the exterior.}\lbl{amalgamate} 
\end{figure}

Now we use the fact that each summand of $\mu$ always contains trivalent
vertices. This follows since the degree~$1$ part is given by linking numbers
which are zero. Therefore, since $\sigma Z(T_i)$ is of the form
$\exp(\text{struts})\exp(\text{rest})$, there is no strut part to $\mu$. (We
need to apply $\sigma$ in order to have an algebra structure.) Therefore, we
have shown that each term of $Z[K\cup L_C;S]$ contains at least $n$
\emph{special} vertices, which by definition are trivalent vertices only adjacent
to other trivalent vertices, or to univalent vertices which lie on internal
Borromean rings components. This also holds for the LMO normalization
$\check{Z}[K\cup L_C;S]$ which only differs by factors of the
Kontsevich integral of the unknot $\nu$.

When we apply $\sigma$ to $\check{Z}[K\cup L_C;S]$, the number of special
vertices can only increase. This can be seen by examining the definition of
$\sigma$ \cite{bn}, which is iteratively defined by operations which involve
removing the skeleton and tacking things on to the created univalent vertices.

Let $S(L_C)$ denote the surgery along $L_C$. Then the Aarhus integral
is defined as
$$
A[(S(L_C),K);S]:=\left(\int^{FG}\sigma\check{Z}[K\cup L_C;S]\right)\cdot
(\text{Kirby I-move normalizations})
$$
The formal Gaussian integration is with respect to the negative inverse of the
linking matrix of the link $L_C$, and this linking matrix is the same in each
summand. Notice that the linking matrix and its inverse are of the forms
$$
\Lambda=\left(\begin{array}{cc}0&I\\I&A\end{array}\right),
\Lambda^{-1}=\left(\begin{array}{cc}-A&I\\I&0\end{array}\right),
$$
where the first row and column refers to B-components and the second to
L-components. We claim that special vertices descend to internal vertices after
applying $\int^{FG}$. If not, then there is a special vertex adjacent to a univalent
vertex labeled by a B component,$b$, and a strut with one endpoint labeled by
$K$ and the other labeled by some component $x$, which are glued together
along a strut labeled by $b$ and $x$, with coefficient coming from
$-\Lambda^{-1}$. By consideration of linking number, the only $K$-$x$ struts
are when $x$ is an L component that links the knot, and hence is not part of an
internal Y. Therefore the $b-x$ strut is between a B component and an L
component in different Y's, which therefore don't link. Thus, by examining
$\Lambda^{-1}$ we see that the weight is zero.

Therefore we have argued that all summands of $A[(S(L_C),K);S]$ have at
least $n$ internal vertices. Multiplying by the normalizations from the Kirby
I-move can only increase this number.   

Notice that for any nonempty $s\subset S$ we have $(S(L_C),K)_s =(S^3,K)$.
Therefore, we have $A[S(L_C),K);S]=A((S^3,K)_C)-A(S^3,K)$. The right hand
side of this last equation is $\sigma Z(K_C)-\sigma Z(K)$, which
has Euler degree~$\geq n$, since we argued that the left-hand side of the
equation has that property. Notice that the wheeling isomorphism can never
decrease Euler degree, since it involves attaching wheels to diagrams. Thus
$\partial_\Omega \sigma Z(K_C)-\partial_\Omega \sigma Z(K)$ is of Euler degree~$\geq
n$.

Finally we must take the logarithm.
Write $a=\partial_\Omega\sigma Z(K_C) - 1$ and $b=\partial_\Omega\sigma Z(K) -
1$. We are interested in 
$$
\log(a+1)-\log(b+1) = \sum_{k=1}^{\infty}(-1)^{k+1}\frac{a^k-b^k}{k}.
$$
 Notice that this is divisible by $a-b=(a+1)-(b+1)$ which we calculated was of Euler
degree~$\geq n$. Since the Euler degree adds under disjoint union (i.e.\
multiplication) it follows that the whole expression is of Euler degree~$\geq n$.
\end{proof} 

\begin{definition}\lbl{def:Z}
Let $Z^g_k :\K \to \B^g_k$ be the grope degree~$k$ part of 
$$
\log_{\hat\B}(\partial_\Omega  \sigma ) Z = (\partial_\Omega  \sigma
)\log_{\hat\A}Z.
$$
\end{definition}

 \begin{corollary} \lbl{gropetypen}
$Z^g_k$ vanishes on $G_{k+1}$.
\end{corollary}

\begin{proof}
Let $K$ be a knot and $C$ be a simple clasper of grope degree~$(k+1)$. 
We need to show that 
$$
z:= \log_{\hat\B}(\partial_\Omega  \sigma ) Z(K)-\log_{\hat\B}(\partial_\Omega 
\sigma ) Z(K_C)\in \widehat\B
$$
has no terms of
grope degree~$\leq k$. Write $k+1=v+b_1=v+\frac{e}{2}+1$ in terms of the Vassiliev
degree and the first Betti number. Then the Euler-degree of $C$ is $2(b_1-1)$
implying by the previous proposition that $z$ starts with terms of  Euler-degree~$2(b_1-1)$. Similarly, by the usual properties of the Kontsevich integral, we know
that $z$ starts with terms of Vassiliev degree~$v$ (this also covers the
case $b_1=0$). Hence our claim follows.
\end{proof}

\begin{lemma}\lbl{splitting}
$Z^g_k\circ \Phi_k^g = Id$.
\end{lemma}

\begin{proof}
Similarly to the above proof, let $U$ be the unknot and 
$C$ be a simple clasper of grope degree~$k$. 
Now we need to show that the $\Z$-linear combination
$$
z:= \log Z(U_C)\in \widehat\B
$$
of diagrams $D_i$ starts (in the grope filtration) with the diagram $D$
underlying the clasper $C$. Writing $k=v+b_1=v(D)+b_1(D)$, we conclude as
in the above argument that
$$
v(D_i)\geq v \text{ and } b_1(D_i)\geq b_1.
$$ 
This implies as before that $g(D_i)\geq g(D)=k$ but also that the grope
degree part of
$z$ consists of the linear combination of those $D_i$ for which $v(D_i)=v$ and
$b_1(D_i)=b_1$. By the usual universality of the Kontsevich integral and
Lemma~\ref{difference}, the first property alone shows that exactly one $D_i=D$
with coefficient
$+1$.
\end{proof}

\begin{corollary}\lbl{bisomorphism}
$Z^g_k$ induces an isomorphism
$$
\G_k\otimes\Q\cong \B^g_k\otimes \Q
$$
\end{corollary}

\begin{proof}
The fact that the map is well-defined is the content of 
Corollary~\ref{gropetypen}. By Lemma~\ref{Phig} the map $\Phi^g_k$ is an
epimorphism. Now Lemma~\ref{splitting} implies that it is injective modulo
torsion, and hence a rational isomorphism. Therefore its rational inverse
$\log Z^g_k$ is also an isomorphism.
\end{proof}

This result clearly implies Theorem~\ref{B} from the introduction.

\subsection{4-dimensional grope cobordism: Grope concordance}

\begin{proposition} \lbl{4D}
For each $k\geq 3$, two knots are class $k$ grope concordant if
and only if their Arf invariants agree.
\end{proposition}

\begin{proof}
It was shown in \cite{cot} that a grope concordance of class $\geq 3$
preserves the Arf invariant. So pick $k\geq 3$ and suppose that $K_1$ and
$K_2$ have the same Arf invariant. We shall construct a grope
concordance of class $k$ as follows: 

Using Lemma \ref{looplem} one can perform clasper
surgeries on claspers {\em with loops} until the two knots share Vassiliev
invariants up to order $k$. 
Hence by \cite[Thm.4]{ct}, the two knots are concordant to knots
$K_1^\prime$ and $K_2^\prime$ that share Vassiliev invariants up to order $k$.
 By \cite[Thm.1]{ct} $K^\prime_1$ and $K^\prime_2$ are then related by
a (capped) grope cobordism of class $k$ in 3-space. This grope cobordism can be
glued to the two concordances to obtain a grope
concordance of class $k$ between $K_1$ and $K_2$.
\end{proof}

\subsection{Clasper moves and the grope degree}\lbl{clgrdeg}
It is the purpose of this section to prove some lemmas on the behavior of
claspers with respect to the grope degree. 

\begin{lemma}\lbl{lemmastar} 
Let $C$ be a rooted tree clasper of degree~$c$ on a knot $K$.
\begin{itemize}
\item[(a)]  Suppose two leaves of $C$ hit $K$ as on the left in Figure
\ref{switchleaves}. Let
$C^\prime$ be obtained from $C$ by interchanging the order of the leaves as
on the right of Figure \ref{switchleaves}. Then $K_C = K_{C^\prime} \mod
G_{c+1}$.
\item[(b)] Suppose $C^\prime$ is a rooted clasper obtained from $C$ by
 by homotoping one of the edges.
 Then $K_C = K_{C^\prime} \mod
G_{c+1}$. Indeed, when $C$ is a tree clasper, $K_C=K_{C^\prime} \mod \T$,
where
$\T$ is formed from the tree type of $C$ by adding a hair to the edge that
is homotoped. 
\item[(c)] Suppose $L$ is a leaf that bounds a disk, and that the leaf has
trivial linking number with $K$. Let
$C^\prime$ be the clasper which has these intersections pushed off of the
disk. Then $K_C = K_{C^\prime} \mod
G_{c+1}$.
\end{itemize}
\end{lemma}

\begin{figure}[ht] 
\begin{center}
\epsfig{file=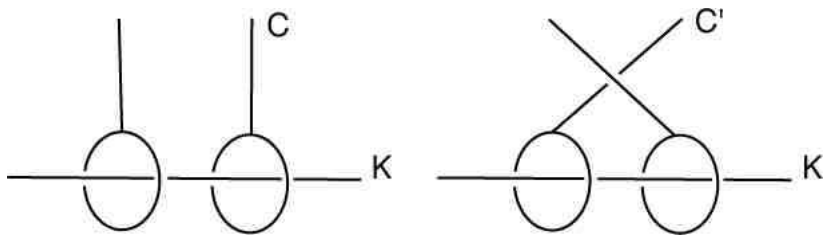,%
          height = 3cm}
\end{center}
\caption{Interchanging the order of leaves.}\lbl{switchleaves}
\end{figure}

\begin{proof}
\cite[Prop.4.4]{h2} proves (a) when the two leaves belong to different
simple claspers.  The hypothesis of  ``simple'' is not used in the proof.
(His Figure 29 is still valid, and the zip construction proceeds without a
hitch.) In order to use this fact, insert some Hopf-linked tips into some
edges of $C$, in order to break $C$ into two claspers, each containing
exactly one of the leaves to be interchanged. This proves (a).

To prove part (b), let $E$ be a degree $1$ clasper with one leaf $L$
linking an edge of $C$ as a meridian, and the other leaf embedded
arbitrarily, so that surgery on $E$ realizes the homotopy of the edge.
Insert a Hopf-pair in the edge of $C$. There are two cases: either this
disconnects $C$ into two claspers $C^\prime$ and $C^{\prime\prime}$ or it
doesn't. Now
$L$ bounds a disk that hits one of the Hopf-pair in two points. Add a tube
to get rid of the intersections. The resulting surface bounding $L$ has a
symplectic basis bounding disks each of which hits one of the Hopf-pair in
one point. In the disconnected case, in the
complement
of $K_C$, these curves therefore bound gropes of the same tree type as
$C^\prime$ and $C^{\prime\prime}$ respectively. Therefore surgery on $E$
is the same as surgery on a clasper formed by gluing the $C^1$ and $C^2$
trees onto the tips of a ``Y". This is exactly the tree type $\T$. In the
connected case, just use one of the symplectic basis elements. 

To prove part (c), consider two intersections of $K$ with the disk of
opposite sign.  Let
$\widetilde{K}$ in
$S^3\backslash C$ be a parallel to $K$ with the two intersection points
pushed off of the disk. 
 Now Figure
\ref{grasper} shows a genus one surface cobounding $K$ and $\widetilde{K}$ in
$S^3_C$. This cobounding surface has a cap which is pierced by the leaf
once. So by Theorem 11 of \cite{ct}, $K$ and $\widetilde{K}$ cobound a class
$c+1$ grope in $S^3_C$ , which says that $K_C$ and $K_C^\prime$ cobound a
class $c+1$ grope in $S^3$. Iterate this procedure until all intersections
are removed.
\begin{figure}[ht] 
\begin{center}
\epsfig{file=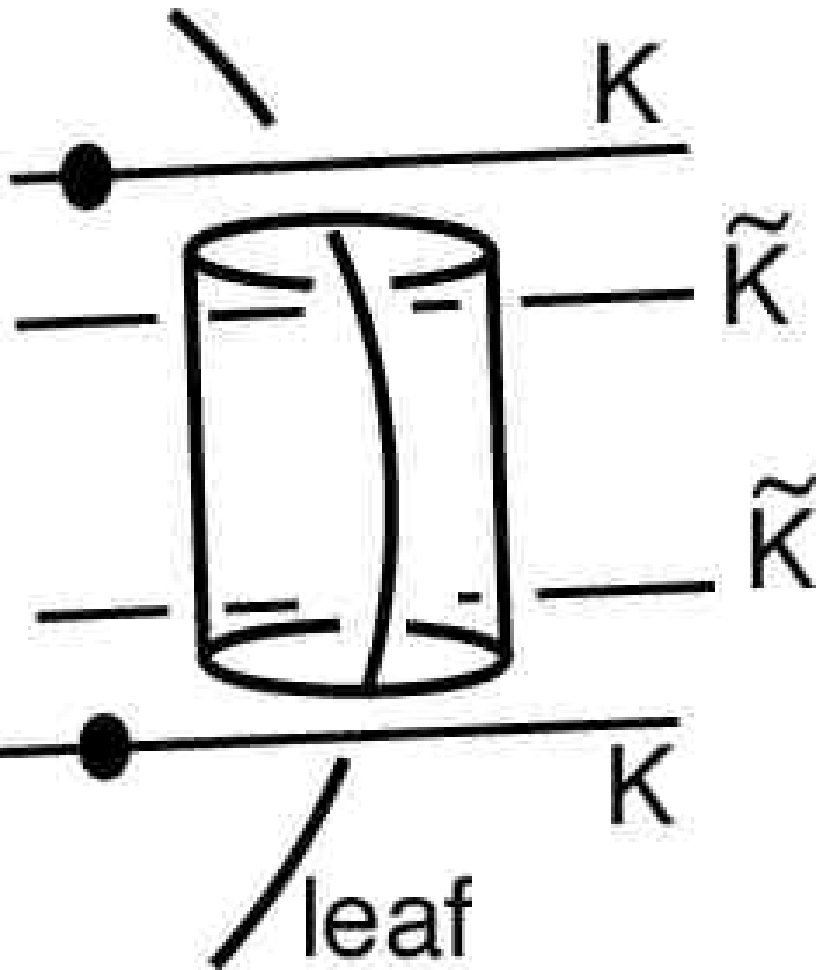,%
          height = 4cm}
\end{center}
\caption{The proof of \ref{lemmastar}(c).}\lbl{grasper}
\end{figure}
\end{proof}

We record here a theorem from \cite{cst}. It is stated in \cite{ggp} for the
case of degree~3 claspers, but the details are ``left to an interested reader."
This theorem is also known to Habiro.
\begin{theorem}[\cite{cst}]\lbl{topihx}
Suppose three claspers $C_i$ represent the three terms in an IHX relation. 
Given an embedding of $C_1$ into a $3$-manifold,there are embeddings of $C_2$ and
$C_3$ inside a regular neighborhood of $C_1$, such that the leaves are parallel
copies of the leaves of $C_1$, and the edges avoid any caps that $C_1$ may
have. Moreover, surgery on $C_1\cup C_2\cup C_3$ is diffeomorphic (rel
boundary) to doing no surgery at all. 
\end{theorem}

\begin{lemma}\lbl{lemmastar2} 
Let $U$ be the unknot.
\begin{itemize}
\item[(a)] Suppose that $C$ and $C^\prime$ are two simple claspers of grope degree~$c$ on $U$ which differ only by a half twist along any edge. Then
$U_C\#U_{C^\prime} \equiv U \mod G_{c+1}$.
\item[(b)] $G_c /G_{c+1}$ is generated under connected sum by knots
which are simple clasper surgeries of grope degree~$c$ on $U$.
\item[(c)] Suppose three claspers $C_i$ of grope degree~$c$ on $U$
 differ according to the IHX relation, see Figure
\ref{1tihx}. Then
$$
U_{C_1}+U_{C_2}+U_{C_3} \in G_{c+1}.
$$
\item[(d)] Suppose two claspers $C^\prime$, $C^{\prime\prime}$ of grope
degree
$c$, on a knot
$K$, differ by a full twist along an edge. Then
$K_{C^\prime}=K_{C^{\prime\prime}} \mod G_{c+1}$.
\end{itemize}
\end{lemma}

\begin{proof}
Part (a):

\noindent  First, insert Hopf-linked pairs of tips to make $C,C^\prime$ trees.
Use Lemma~\ref{preinverse} to find $\widetilde{C}$, such that
$U_{C\cup\widetilde{C}}=U$. 
 We need to disentangle $\widetilde{C}$. Those leaves that were meridians to
$U$ on $C$, are still meridians on $C^\prime$. For every tree clasper $D$
of degree exceeding $1$ on a knot $K$, a meridian to each (non-root) leaf
and edge links
$K_D$ algebraically trivially. Hence the pushed off Hopf-linked tips on
$\widetilde{C}$ link $U_C$ algebraically trivially, and so all intersections
can be pushed out by Lemma
\ref{lemmastar}(c). The pushed off Hopf- linked tips of $\widetilde{C}$, are now
Hopf-linked in the same way as $C$. Denote the new clasper
$\widetilde{C}^\prime$.
 Let
$B$ be a ball meeting $U_C$ in a standard unknotted arc away from $C$. Slide
the leaves of
${\widetilde{C}}^\prime$ into $B$. By Lemma~\ref{lemmastar}(b), we can pull
${\widetilde{C}}^\prime$ into $B$ modulo $G_{c+1}$. Let this new clasper be
called $C^\prime$. We have just demonstrated that
$$
U \equiv (U_C)_{\widetilde{C}} \equiv U_C\#U_{C^\prime} \mod G_{c+1}
$$
and that $C^\prime$ is of the required form.

\noindent Part (b):
By Theorem 2 of \cite{ct}, if $K\in G_c$, then there are knots $K_i$ and
simple claspers $C_i$ of grope degree~$\geq c$, such that $K_0=U, K_i =
\left(K_{i-1}\right)_{C_i}, K_N=K$. 
Modulo $G_{c+1}$, we can discard all $C_i$ except those of degree~$c$.
As in part a, slide the leaves of $C_N$
into a ball on $K_{N-1}$, and then pull $C_N$ into the ball by Lemma
\ref{lemmastar}(b). Call the new clasper
$C_N^\prime$. Then $K_N \equiv K_{N-1}\#U_{C_N^\prime} \mod G_{c+1}$.
Inductively, we are done.

\noindent Part (c): 
By Theorem~\ref{topihx}, we can find three tree claspers
differing by the IHX relation inside a regular neighborhood of each other, such that
surgery on all three is null-isotopic.
 As noted
previously, we can pull them apart modulo higher grope degree. After pulling them
apart Hopf pairs of leaves can be blown down into edges. 

\noindent Part (d):
There is a clasper $\widetilde{C}$ that differs from both $C^{\prime}$ and
$C^{\prime\prime}$ by a single half-twist. Then  the proof of
part (a) implies that $K_{C^\prime} = (U_{\widetilde{C}})^{-1}\mod G_{c+1}$,
and similarly $K_{C^{\prime\prime}} = (U_{\widetilde{C}})^{-1}\mod G_{c+1}$. 
\end{proof}

\section{Low degree calculations} \lbl{sec:low}

All the $4$-dimensional results in Table~\ref{table} are contained in
\cite{cot}, so we work purely in dimension~$3$.
 We begin by making a general
observation which turns out to be very useful for all our calculations.
\begin{proposition}\lbl{r1}
$$ K_1 \equiv K_2 \in \K/G^c_{k} \text{ implies } K_1\equiv K_2
\in\K/G_{k}\text{ implies } K_1\equiv K_2 \in\K/G^c_{
\lfloor{\frac{k}{2}\rfloor+1}}
$$
\end{proposition}

\noindent{Remark: The second implication is a slight improvement over
\cite{newconant}, although its proof depends on Theorem 2 of \cite{ct}, which
in turn depends on
\cite{newconant}.} 

\begin{proof}The first implication is obvious: a capped grope is also an uncapped grope.
The second implication arises as follows. A simple clasper of grope degree~$k$ has minimal Vassiliev
degree when the first Betti number is maximized. Suppose $k=2t+1$. Then the number of
leaves of the clasper when the edges are cut to make a tree is $2t+2$. Note that the
grope degree is unchanged under performing such cuts. At most $2t$ of these edges
could have been paired together to form
$t$ loops. The Vassiliev degree is then
$k-t =\lfloor{\frac{k}{2}}\rfloor +1$. 

If $k=2t$, then there are $2t+1$ leaves of the associated tree clasper,
$2t$ of which can be paired to make a tree. In that case, we would have a
clasper of degree~$t$, with a single leaf hitting the knot. Since the
corresponding Feynman diagram is trivial modulo STU, the map $\Phi$
indicates that it must be trivial modulo $G^c_{t+1}$.   
\end{proof}

Moreover, the preceding argument proves that the only type
 $\lfloor k/2 \rfloor +1$ invariant values that can be attained by a
$G_k$-trivial knot are  those corresponding to linear combinations of
connected Feynman diagrams with $\lfloor  (k-1)/2 \rfloor$ loops, using Lemma
\ref{difference}.

\subsection{The groups $\K/G_k$ for $k\leq 5$}

\begin{theorem}\lbl{calc} We have the following calculations:
\begin{align*}
\K/G_1 &= \{ 0\} &&\K/G_2 = \{ 0\}\\
\K/G_3 &\cong \Z/2 (c_2) &&\K/G_4 \cong \Z (c_2)\\
\K/G_5 &\cong \Z(c_2)\oplus \Z/2(c_3)
\end{align*}
\end{theorem}

Here $c_2$ and $c_3$ denote some choice of the degree~$2$ and $3$
Vassiliev invariants.
The proof uses the following well-known calculations of the indecomposable
elements
$\A^I$, see for instance \cite{g} or \cite{bn}. The last statement is
due to \cite{ng}.

\begin{lemma}\lbl{looplem}
$\A^I_2\cong \Z$, with generator \epsfig{file=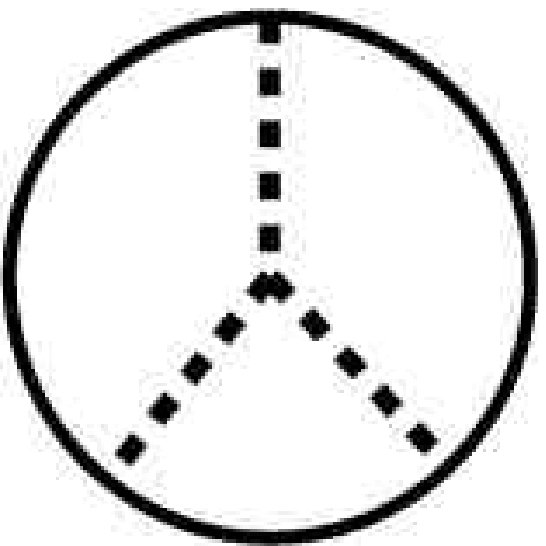, height=.5cm}
 and $\A^I_3\cong \Z$ with generator \epsfig{file=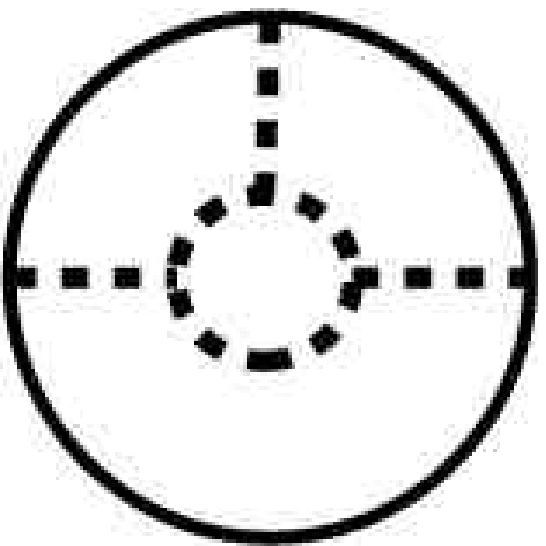,
height=.5cm}. In general, $\A^I_k$ is generated by connected diagrams with at
least $1$ loop for $k\geq 3$.
\end{lemma}

\begin{definition} Let $\A^I_k[m]$ denote $\A^I_k$ modulo diagrams with
$m$-loops.
\end{definition}

\begin{proof}[Proof of \ref{calc}]
Suppose $\A^I_m$ is torsion free for $m\leq k$. 
We have seen that the maps $\Phi_m:\A^I_m\to\G_m^c$ from Lemma~\ref{Phi} are then
isomorphisms for $m\leq k$.
Hence 
$$
\K/G^c_k \cong \oplus_{m < k}\A^I_m
$$
By Proposition~\ref{r1}, $\K/G_k$ is a quotient group of
$\K/G^c_k$  by simple clasper moves of grope degree~$k$. A
simple clasper move represents a diagram in
$\A^I_l$ for some $l$, and by definition the number of loops is $n-l$. Hence, by
Lemma~\ref{difference}, the degree~$l$ part changes according to the
corresponding Feynman diagram. Thus we have a relation of the form
$$(0,\ldots,0,\alpha,*,\ldots,*) \in
\A^I_2\oplus\A^I_3\oplus\cdots\oplus\A^I_l\oplus\cdots\oplus
\A^I_k$$ where $\alpha$ lies in the subspace of $\A^I_l$ generated by
diagrams with $n-l$ loops. These constitute all of the relations, but
notice that we have no control over the $*$'s, so this fact will only be
useful in (very) low degrees.

Now $\K/G_2$ is a quotient of $\K/G^c_2= \{0\}$, and is therefore
trivial. That $\K/G_2$ is trivial also follows from the statement that all
knots cobound a surface with the unknot. For instance, one may take a
punctured Seifert surface.

By the above remarks $\K/G_3$ is a quotient of $\A^I_2$, by the subspace of
diagrams with
one loop. The $1$-loop subspace is generated by the following diagram,
which, as shown, is equal to twice the generator.

\begin{center}
\epsfig{file=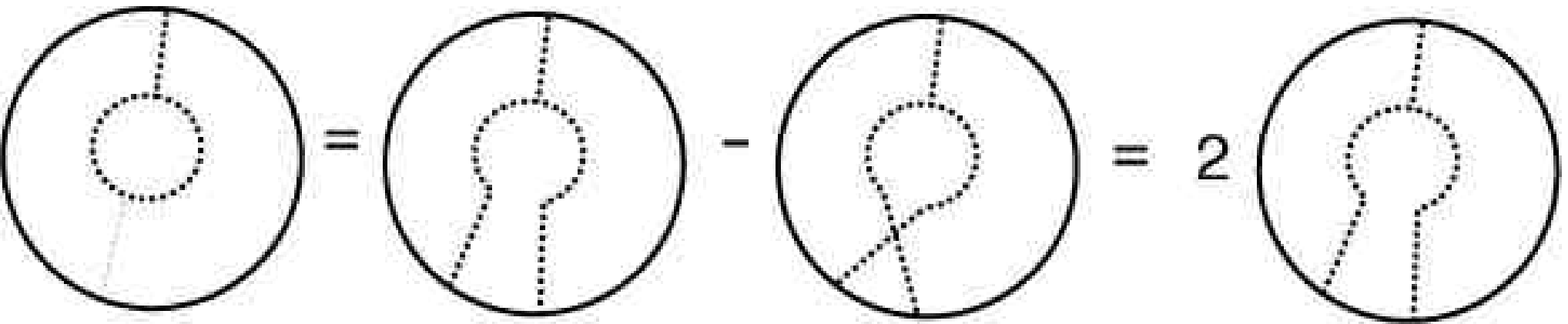,%
        height = 2cm}
\end{center}

Since the Arf invariant is the mod 2 reduction of $c_2$, we have proven what
we needed to.
For the next degree, note that $\K/G_4$ is a quotient of $\A^I_2\oplus \A^I_3$
by relations of the form:

\begin{center}
 (2-loop, $*$)  and  (0, 1-loop)
\end{center}

By Lemma~\ref{looplem}, the relations of the second type kill off
$\A^I_3$, and we are left with 
$A^I_2[2]$. The $2$ loop subspace of $A^I_2$ is generated by a diagram with
one foot on the solid circle, which is trivial modulo $STU$. Hence
$\A^I_2[2]=\A^I_2$.

Now $\K/G_5$ is a quotient of $\A^I_2\oplus\A^I_3\oplus\A^I_4$, by relations
of the form

\begin{center}
 (0 , 2-loop, $*$)  and  (0 ,0, 1-loop)
\end{center}

Again, by \ref{looplem} the
second type of relation kills off $\A^I_4$, and we are left with
$\A^I_2\oplus\A^I_3[2]$.
By \cite{ng}'s arguments the $2$-loop subspace is generated by diagrams of
the form
\epsfig{file=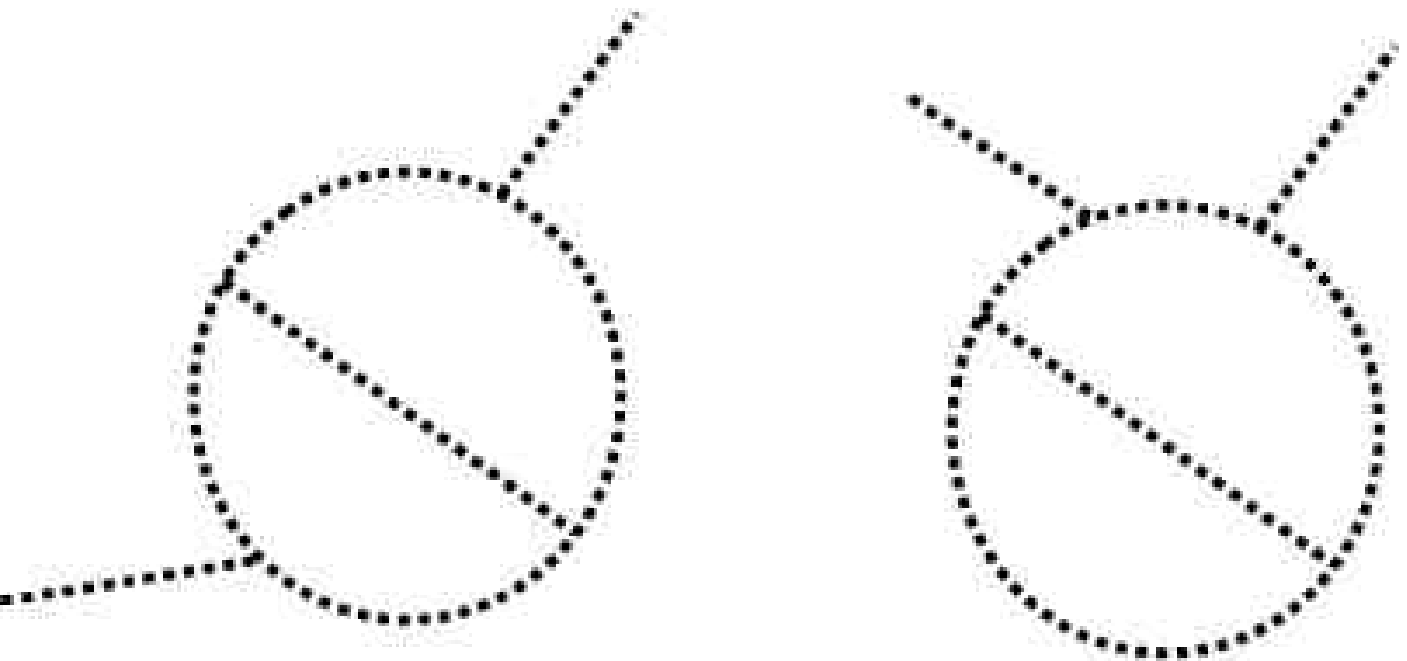,%
        height = .5cm}
 attached to the outer loop by some permutation. For instance
we do not need to separately consider the diagrams
\epsfig{file=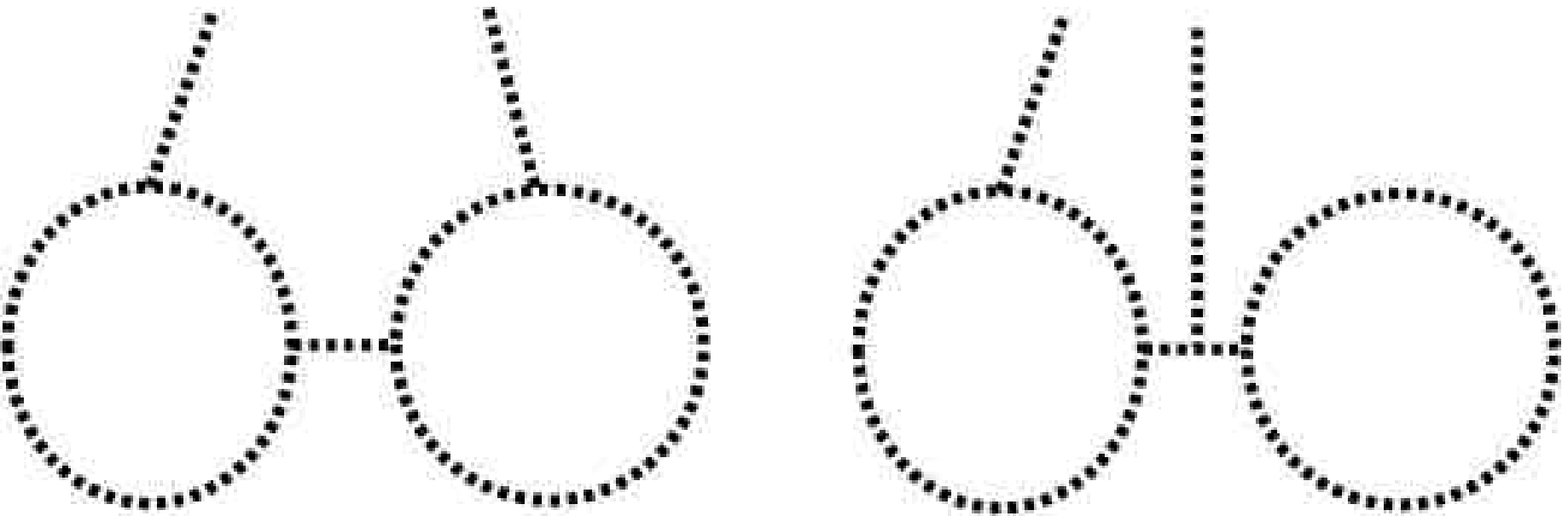,%
        height = .5cm}.
 Up to sign, there is only one diagram of each of the good
types, and they can each be represented as follows.
\begin{center}
\epsfig{file=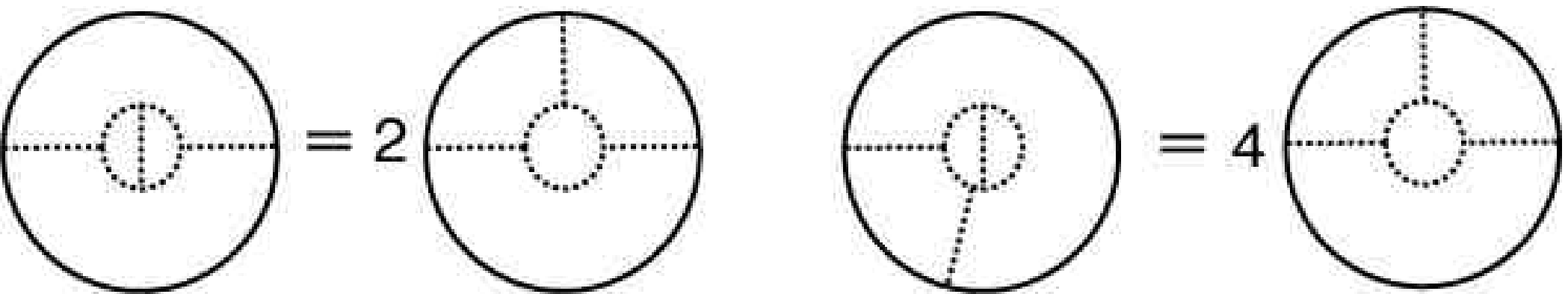,%
        height = 2cm}
\end{center} Since twice the generator is realized, but the generator is
not realized,
$c_3
\text{ mod } 2$ is all that survives.
\end{proof}

The previous calculations are quite suggestive of the following conjecture.

\begin{conjecture} \lbl{integral}
Suppose $\A^I$ is torsion free up to degree~$k$. Then 
$$
\K/G_k \cong \A_{k-1}[1]\oplus\A_{k-2}[2]\oplus \A_{k-3}[3]\oplus
\cdots\oplus\A_2[k-2].
$$ 
\end{conjecture}

We remark that this is true rationally, without the hypothesis, see
Corollary~\ref{bisomorphism}. We end this section with a generalization of one of
our results above to knots in integral homology spheres.

\begin{proposition} 
Two knots in an integral homology sphere have the same
Arf invariant iff they are class $3$ grope cobordant.
\end{proposition}

\begin{proof}
Let $K$ be a knot in an integral homology sphere $M$. It bounds a surface.
By Matveev's result \cite{m}, there is a collection of Vassiliev degree~$2$
claspers which turn $M$ into
$S^3$. We may assume that the claspers are disjoint from the surface since
we can perturb them by isotopy. One can then find inverse claspers in a
regular neighborhood, so that we are in the situation of a knot in $S^3$
together with some claspers disjoint from a Seifert surface, such that
surgery on these claspers takes us to the original pair $(M,K)$. Now by
Lemma \ref{lemmastar}(c), $K$ is class $3$ cobordant to a knot in
$D^3\subset M$. Now we can use the result for $S^3$.
\end{proof}

\subsection{Tree types of class~$4$}

In Theorem 4.2, we analyzed the equivalence relation given by grope
cobordisms of a fixed class, up to class $5$. In general, when one refines
these equivalence relations to be of a fixed tree type, one gets a
different answer. For instance, $S$-equivalence is generated by a specific
tree type of class $5$ (Theorem \ref{symtheorem}). Degree $4$ is the first
place that the rooted tree type is not unique, but we prove in this section
that class $4$ cobordism is generated by either of the two rooted tree
types: 
\epsfig{file=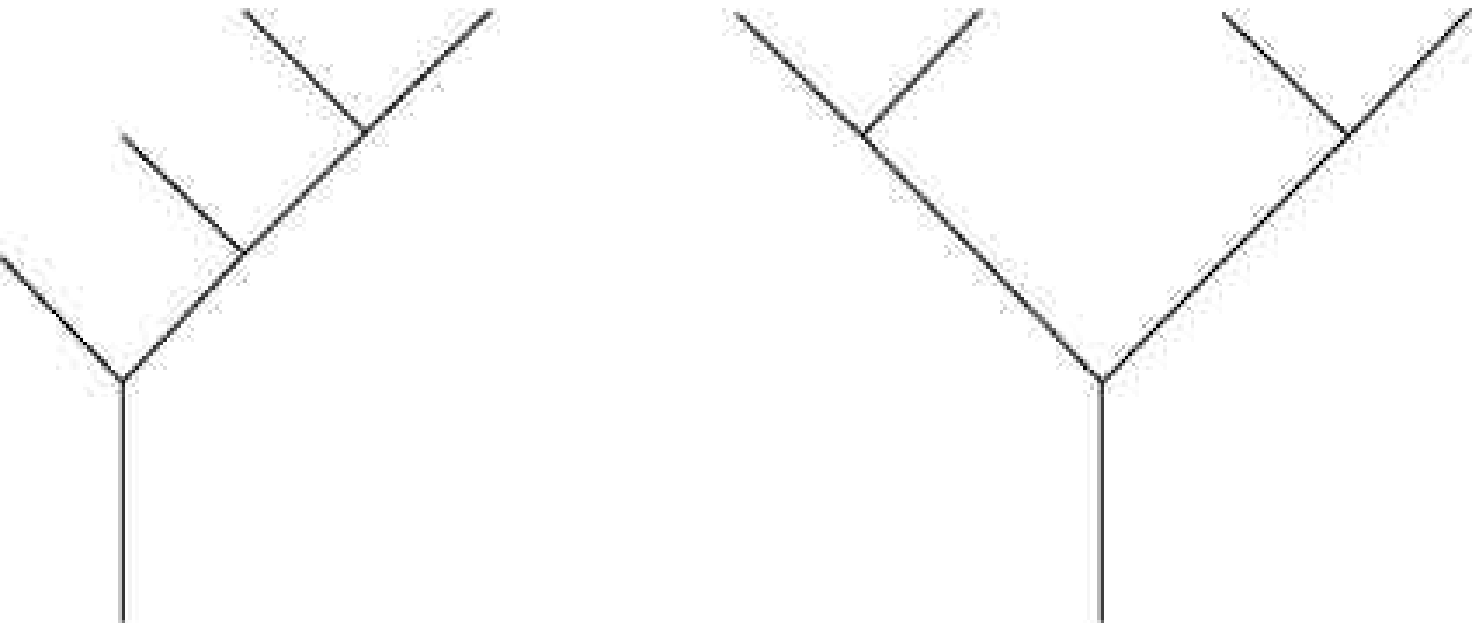, height = .4cm}.

\begin{proposition}\lbl{whoze}
$\K/\epsfig{file=treetypes4half.eps,%
        height = .4cm} \cong
\K/\epsfig{file=treetypes4full.eps,%
        height = .4cm}
 \cong \K/G_4$. That is, grope cobordism of class four is generated
by either of the two rooted tree types.
\end{proposition}

\begin{proof}
Let the two tree types be called $T_1$ and $T_2$. We have the following
commutative diagram.
\begin{equation*}
\begin{CD}
0@>>> Ker @>>> \K/G^c_4 @>>> \K/G^4 @>>> 0\\
@AA=A @AA{onto}A @AA=A @AA{onto}A @AA=A\\
0 @>>> Ker_i @>>> \K/G^c_4 @>>> \K/T_i @>>> 0
\end{CD}
\end{equation*}
As we saw in the proof of Theorem~\ref{calc}, the $Ker$ is generated by a
wheel with three legs attached to the outer circle. After cutting this, the
root can be chosen so that it is either of the two tree types. That implies
the map
$Ker_i\to Ker$ is onto. By a the five lemma the map $\K/T_i\to K/G_4$ is
an isomorphism.
\end{proof}

\subsection{Tree types of class~$5$}

In this section we prove that grope cobordism of class $5$ is generated by
\epsfig{file=half5.eps,%
        height = .4cm} - equivalence, or by 
\epsfig{file=treetypes5anom.eps,%
        height=.4cm} -equivalence. In \cite{ct}, it was proven that ``half
gropes" generate equivalence by (capped) gropes of a given class. Here a
half grope has a tree type representing a right-normed commutator.  This
implies the first isomorphism in the following lemma. 

\begin{lemma}\lbl{l3}
$\K/\epsfig{file=half5.eps,%
              height=.4cm} \cong
\K/G_5
\cong
\K/\epsfig{file=treetypes5anom.eps,%
        height = .4cm}
$
\end{lemma}

\begin{proof}As in the proof of Lemma~\ref{whoze}, it suffices to show that
the kernel of $\K/G^c_5\to\K/G_5$ is generated by both of the tree types.
This kernel is generated by a circle with four legs, something that can be
thought of as either tree type, as well as twice the
generator corresponding to
$c_3$. This is realized by a theta with two legs which can be cut apart to be
either of the two tree types. \end{proof}

\subsection{$S$-equivalence}

Let $\mathcal S_5$ denote the following tree type:
\begin{center}
\epsfig{file=sym5.eps,%
        height = 1cm}
\end{center} This is the simplest tree with an internal vertex.

\begin{theorem}\lbl{symtheorem} Two knots are $S$-equivalent iff they are
$\mathcal S_5$-equivalent.
\end{theorem}

\begin{remark}
In \cite{cot} it is proven that the corresponding move in $4$
dimensions gives Blanchfield forms up to cobordism. Thus the kernel of
going from 3 to 4 dimensions consists of adding the relation
$K+K^!=0$ where $K^!$ is the mirror image and is an inverse in the knot
concordance group.
\end{remark}

\begin{proof}
That $\mathcal S_5$ preserves $S$-equivalence follows by a construction of
Murakami and Ohtsuki \cite[p.6]{mo}, applied to a disk leaf of the $\mathcal
S^{3D}_5$ move. 
For the converse we use a result of Naik and Stanford\cite{ns} (see also
Murakami and Nakanishi \cite{mn}), that the
\emph{doubled delta move} generates $S$-equivalence.
\begin{figure} 
\begin{center}
\epsfig{file=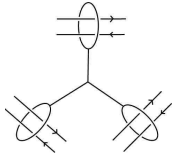,%
        height = 2cm}
\end{center}
\caption{Stanford and Naik's doubled delta move.}\lbl{ddelta}
\end{figure} In fact, \cite{ns} proves the stronger result that the
doubled delta move applied to bands of some Seifert surface generates
$S$-equivalence. If this move is applied to three bands, no two of which
are dual, it is easy to construct an $\s_5$ grope cobounding the knots
before and after the doubled delta move. We construct this grope by
constructing disjoint surfaces bounded by the three leaves of the doubled
delta move. These surfaces are constructed by tubing into the dual band.
See Figure
\ref{tube}.
\begin{figure} 
\begin{center}
\epsfig{file=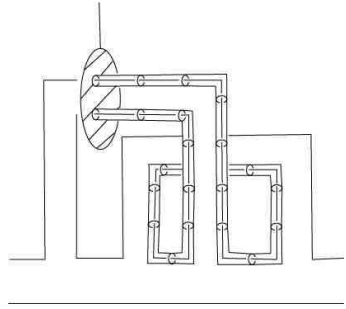,%
        height = 4cm}
\end{center}
\caption{Tubing into the dual band.} \lbl{tube}
\end{figure}
If two or three leaves link the same band, we can construct the surfaces
by nesting the tubes.

 Now that we have these $3$ surfaces, we get an
$\s_5$ grope by the discussion of grope-clasper duality given in \cite{ct}.
A meridian of an innermost tube provides the root.
We are left with the following two cases: 
\begin{itemize}
\item[(a)] Two leaves link one band, and
the other leaf links the dual band. 
\item[(b)] Two leaves link dual bands, and the
third leaf links some other band.
\end{itemize}
Note that in both cases we can construct disjoint surfaces on two of the leaves. 
The argument that $S$-equivalence is generated by the 
$\s_5$ grope is given by diagram \ref{whocares}.
\begin{figure}[ht] 
\begin{center}
\epsfig{file=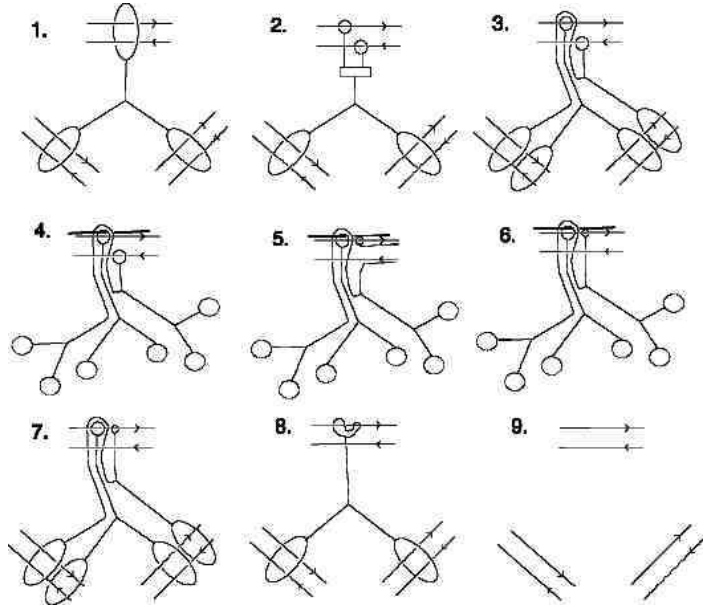,%
        height = 8cm}
\end{center}
\caption{The proof that the $\s_5$ grope generates
$S$-equivalence.}\lbl{whocares}
\end{figure}
We now give the explanation of the figure. First, the two
leaves which bound surfaces are the leaves at the bottom of frame~1. When
constructing these surfaces we can assume that they run parallel to only one of the
strands going through the top leaf. Let this strand be the top strand of the
picture. Going to frame 2 is Habiro's move $8$. Going to frame $3$ follows easily
from Habiro's move 11. In frame $4$ we have added the two surfaces to the tips of
one of the claspers. Rather than drawing the details, we have represented things
schematically. The heavy line at the top of the picture reminds us that pieces of
the new clasper run through here. In frame $5$ we have slid one of the leaves along
the knot until it is opposite the other leaf. Because of the fact that the strands
were oriented oppositely, the leaves face each other as shown. Going to
frame $6$ is where we use the $\mathcal S_5$ move. We pull the edge of the
clasper which is incident to the leaf we just slid into the shown position.
Then by Lemma~\ref{lemmastar}(b), the knots in pictures
$5$ and $6$ cobound an $\s_5$ grope. To get to frame $7$ we forget the fact
that there are surfaces bounding two of the leaves. (Basically the inverse
of
$3\to 4$.) In frame $8$ we do Habiro's move 11 in reverse. But then we have
a clasper with trivial leaf, which does not alter the knot as in frame $9$.
\end{proof}

Note that Corollary~\ref{cor:c3} follows immediately, because by
Theorem~\ref{calc},
$c_3 \mod 2$ is an invariant of all class $5$ grope cobordisms.

This corollary is somewhat surprising. \cite{mo} have proven that the only
rational finite type $S$-equivalence invariants are the coefficients of the
Alexander-Conway polynomial. Recall that $c_3 \mod 2$ is not an
Alexander-Conway coefficient. Recently, Ted Stanford \cite{st} discovered that
$c_3 \mod 2$ can be expressed as a polynomial in the Conway coefficients
$c_2$ and $c_4$. Hence all known finite type invariants of
S-equivalence come from the Alexander polynomial. 

We have classified the behavior of all capped and uncapped class $5$ trees
but the following.

\begin{problem} What is $\K/(\s_5)^c$?
\end{problem}

Clearly a $(\s_5)^c$ move must preserve $S$-equivalence and also type $4$
Vassiliev invariants. Conjecturally this completely characterizes the move.

We close this section with a natural conjecture based on our low-degree
calculations.

\begin{conjecture}
$\K/\T$ only depends on the unrooted tree type of $\T$.
\end{conjecture}

\section{Miscellaneous results}\lbl{sec:misc}

\subsection{Null filtration of knots which bound a grope} \lbl{sec:null}

In the first author's Ph.D. thesis \cite{c} (see also \cite{newconant}), the
question of a knot
\emph{bounding} a grope is considered. This is much stronger than
cobounding a grope with the unknot. In particular, if a knot bounds a grope
of class at least three, the knot has trivial Alexander polynomial, whereas
this is certainly not the case for the cobounding situation. The central
result in \cite{c} is the following:

\begin{theorem} If a knot $K$ bounds an embedded grope of class
$k$ into $S^3$, then Vassiliev invariants up to degree $\lceil k/2\rceil$
vanish on $K$.
\end{theorem}

In fact the bounding of a grope is an extremely restrictive condition. This can be
very well expressed in terms of the {\em null filtration} of \cite{gr}. It is
obtained by the usual alternating sum formalism by declaring a null clasper surgery
on a Y to have degree one. Here the word ``null'' expresses the condition that the
leaves of the clasper must have trivial linking numbers with the knot. It follows
that to be null equivalent to the unknot (the case $k=3$ below) is the same as
having trivial Alexander polynomial, at least in a homology sphere.

\begin{proposition} \lbl{bounding}
If a knot $K$ bounds an embedded
grope of class $k$ in a $3$-manifold $M$, then the pair $(M,K)$ is
$(k-3)$-null equivalent to the unknot in $M$.
\end{proposition}

\begin{proof}
A knot bounding a grope of class $k$ can be obtained from the unknot $U$ by surgery
on a rooted tree clasper of Vassiliev degree~$k$, where $U$ is a 
meridian to the root. In particular, the other leaves do not link $U$.

Break the
clasper into a union of
$(k-1)$ Y's. Surger the knot along the Y which contains the root. This
leads to a union of $(k-2)$ Y's, which the surgered knot links trivially.
Now consider the alternating sum, surgering over all subsets of the Y's.
It is easy to see that surgery on any proper subset does not change the
knot, because there will be leaves which bound embedded disks. Hence the
alternating sum reduces to $(M,U)-(M,K)$ and since we did $(k-2)$ Y-surgeries,
these two knots are $(k-3)$-null equivalent.
\end{proof}

\subsection{Grope cobordism and orientation reversal} 

\begin{proposition} \lbl{2-torsion} 
Let $\rho$ be the map reversing a knot's
orientation. Then for every knot $K$ in the $k$-th term $G_k$ of the grope
filtration of $\K$, one has 
$$
K\equiv (-1)^k\rho(K)\mod G_{k+1}.
$$
\end{proposition}

 One can filter the primitive Feynman diagrams by grope degree, and in the
associated graded group, it is straightforward to show that
$D=(-1)^{|D|}\rho(D)$.  Conjecturally the graded pieces of this group are
isomorphic to $\G_k$, in which case we'd be done. However, we can still
mimic the Feynman diagram computation geometrically. 

\begin{figure}[ht] 
\begin{center}
\epsfig{file=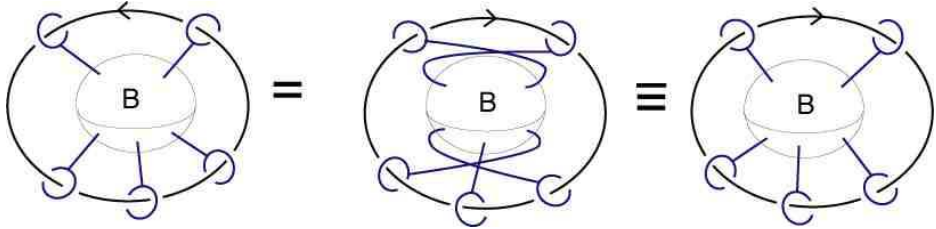,%
        height = 3cm}
\end{center}
\caption{The proof of Proposition~\ref{2-torsion}.}\lbl{fig:flip}
\end{figure} 

\begin{proof}[Proof of Proposition~\ref{2-torsion}]
By Lemma~\ref{lemmastar2}(b) and Lemma \ref{lemmastar}(b), it suffices to show the
result for $K=U_C$, where $C$ is a simple clasper of
degree~$k$ which hits the unknot $U$ as on the left of Figure~\ref{fig:flip}. Here
$B$ is a ball which contains most of the clasper, excluding the leaves. 
 One rotates the outer unknot about the vertical axis lying in
the page, while keeping the ball $B$ fixed. By Lemma \ref{lemmastar}, we
can reorder the leaves as on the right of the picture. Each edge incident
to a leaf has picked up a half twist. Note that the parity of the number of
leaves of a simple clasper matches the parity of the grope degree~$k$.
By Lemma~\ref{lemmastar2}(a), we have shown that $U_C=(-1)^k\rho(U_C) \in
G_k/G_{k+1}$. 
\end{proof}

\subsection{Simple clasper surgeries and the grope degree} 
In this section, we restrict attention to knots and claspers in 3-space.
Consider the statement, 

\noindent ``A simple clasper surgery of degree~$k$ may be realized
by a sequence of simple clasper surgeries of degree~$(k-1)$." 

When the degree~$k$ is
the Vassiliev degree, this follows from Habiro's work: By his main theorem any
simple clasper surgery of Vassiliev degree~$k$ may be realized by a sequence of
simple {\em tree} clasper surgeries of Vassiliev degree~$k$. 
 For a simple tree claspers, one breaks the
clasper into a union of a  degree~$(k-1)$ tree clasper and a Y, with two
leaves linking as a Hopf pair (exactly one of the leaves belongs to the Y). Surgering
the knot along the Y, we get a capped tree clasper of Vassiliev degree~$(k-1)$,
which can then be refined via the zip construction into a sequence of simple tree
clasper surgeries of Vassiliev degree~$(k-1)$. 

It is the purpose of this section to demonstrate that the above statement holds
also for the grope degree.

\begin{theorem}\lbl{imply} 
A simple clasper surgery of grope degree~$k$
may be realized by a sequence of simple clasper surgeries of grope degree~$(k-1)$.
\end{theorem}
We make the following preliminary definition for the purpose of this
section.
\begin{definition} A clasper is said to be \emph{admissible} if it is
simple and the graph type is obtained from a connected trivalent graph with
no separating edges by adding a positive number of legs to the edges.
We also call the trivalent graph admissible.
\end{definition}

Let $\mathcal S$ denote the set of finite sequences of connected
uni-trivalent  graphs with at least one univalent vertex. Consider the
partial order generated by the relation
$$(G_1,\ldots, G_i,\ldots, G_N) < (G_1,\ldots, G_i^1, G_i^2,
G_{i+1},\ldots, G_N),$$ where $G_i^j$ are the two other terms in an IHX
relation involving $G_i$. This partial order gives a convenient language to
state the following lemma. 

\begin{lemma}\lbl{simple} For all $s\in \mathcal S$, there is an $s_0\in
\mathcal S$ such that $s<s_0$ and $s_0= (G_i)$ is a sequence of admissible
graphs
$G_i$.
\end{lemma}

\begin{proof}
Straightforward.
\end{proof}

The point is that these elementary relations can be realized topologically by
embedded claspers, as the next proposition makes clear.

\begin{proposition}\lbl{realize}
Let $G_1,G_2,G_3$ be uni-trivalent graphs related by an IHX
relation, and let $C_1$ be a simple clasper of type $G_1$, embedded in the
complement of a knot $K$. Then $C_1$ may be realized by a sequence of claspers one of
which is of type $G_2$, one of which is of type
$G_3$, and the rest of which have increased grope degree and (at least) the
same number of simple caps.
\end{proposition}

\begin{proof}  First, we must convert some edges of $C_1$ to Hopf
pairs to obtain a tree clasper $\hat{C}_1$ of tree type $\hat{G}_1$.
This gives rise to induced tree types $\hat{G}_2, \hat{G}_3$. There is an
inverse $\overline{C}_1$ to the clasper $\hat{C}_1$ inside a regular
neighborhood of $\hat{C}_1$. The leaves of the inverse are parallels of the
original, but the edges may wander around the regular neighborhood in a
complicated way. The edges do however avoid any caps that $\hat{C}_1$ may
have.(See Lemma~\ref{preinverse}.) If $K_1$ is the knot after surgery on
$C_1$ (equivalently on
$\hat{C}_1$) then $\overline{C}_1$ sits on the knot $K_1$, and surgery on
it produces the original knot. The knot $K_1$ will wander through the Hopf
pairs of $\overline{C}_1$, but will link these leaves trivially. By Lemma
\ref{lemmastar}(c), the surgery $\overline{C}_1$ can be realized modulo higher grope
degree, by surgery along a clasper $\widetilde{C}_1$ obtained from $\overline{C}_1$ by
pushing the knot out of the Hopf pairs. 

Now we use the topological IHX relation, Theorem~\ref{topihx}.
 There is a union of two tree claspers $\hat{C}_2$ and $\hat{C}_3$ in a
regular neighborhood of $\widetilde{C}_1$ which are of type $\hat{G}_2$ and
$\hat{G}_3$ respectively, such that surgery on $\widetilde{C}_1\cup
\hat{C}_2\cup\hat{C}_3$ is null isotopic.  The leaves of
$\hat{C}_2$ and $\hat{C}_3$ are parallels of the corresponding leaves of
$\widetilde{C}_1$. The edges of $\hat{C}_2$ and $\hat{C}_3$ may run through the
regular neighborhood of
$C_1$, but avoid any caps that $C_1$ may have. 

Where there were Hopf pairs of leaves on $C_1$, all the leaves of the three
claspers link. However, these leaves of $\hat{C}_3$ have trivial
linking number with $(K_1)_{\widetilde{C}_1\cup \hat{C}_2}$. Thus, by Lemma
\ref{lemmastar}(c),
$K_1 =((K_1)_{\widetilde{C}_1\cup \hat{C}_2})_{\hat{C}_3}$ is equivalent modulo
higher grope degree to $((K_1)_{\widetilde{C}_1\cup
\hat{C}_2})_{\widetilde{C}_3}$ where $\widetilde{C}_3$ is a clasper obtained from
$\hat{C}_3$ by pushing strands of the knot out of the Hopf pairs of leaves.
Similarly, modulo higher grope degree $\hat{C}_2$ can be realized by a
clasper
$\widetilde{C}_2$, where the knot has been pushed out of the Hopf pairs. Thus
$K_1$ is equivalent modulo higher grope degree to 
$(((K_1)_{\widetilde{C}_1})_{\widetilde{C}_2})_{\widetilde{C}_3}$. However
$(K_1)_{\widetilde{C}_1}$ is equivalent modulo higher grope degree to $K$.
Hence $(K_{\widetilde{C}_2})_{\widetilde{C}_3}$ is equivalent modulo higher grope
degree to $K_1$. The claspers $\widetilde{C}_2$, $\widetilde{C}_3$ have clean Hopf
pairs, and when these are converted to edges, the resulting simple claspers
have graph type $G_2$ and $G_3$ respectively.\end{proof}

Recall that a cap is an embedded disk bounding a leaf of a clasper $C$, with
interior disjoint from $C$. If there are several caps, they are assumed to be
embedded {\em disjointly}. A cap is simple (with respect to a knot $K$) if it has a
single intersection with $K$.
\begin{definitions} Let $C$ be a clasper having some caps.
\begin{itemize} 
\item[(a)] Let $c_1(C)$ be the number of simple caps.
\item[(b)] Let $g(C)$ be the grope degree of the clasper $C$.  After breaking some
edges into Hopf pairs of leaves to make $C$ a tree clasper, $g(C)$ is the number of
leaves minus one.
\end{itemize}
\end{definitions}

\begin{lemma}\lbl{lemmagraphmoves}
Consider a surgery on a clasper $C$.
\begin{itemize} 
\item[(a)] It may be realized by a sequence of
surgeries on claspers $C_i$ which are of the following two possible forms.
\begin{itemize}
\item[(i)] $g(C_i) = g(C)$, $C_i$ is admissible, $c_1(C_i)\geq c_1(C)$
\item[(ii)] $g(C_i) = g(C)+1$, $c_1(C_i)\geq c_1(C)$.
\end{itemize}

\item[(b)] If $C$ is admissible then the clasper surgery may be realized by a
sequence of surgeries on claspers $C_i$ of the following two possible forms.
\begin{itemize}
\item[(i)] $g(C_i) = g(C)-1$, $c_1(C_i)> c_1(C)$
\item[(ii)] $g(C_i) = g(C)$, $c_1(C_i) > c_1(C)$
\end{itemize}
\end{itemize}
\end{lemma}

\begin{proof}Part (a):

\noindent By Theorem 20 of \cite{ct}, we may realize $C$ by a sequence of {\em
simple} clasper surgeries of the same grope degree and at least the same number of
simple caps together with some surgeries of higher grope degree which have at least
the same number of simple caps. (One must check that \cite{ct} Lemma 19
does not decrease the number of simple caps. Some caps may be destroyed
when adding the nested tubes at one stage of the proof. However, instead of
adding the nested tubes, push the disk over the cap as we did for the root
leaf.)

 Hence it suffices to prove (a) for simple claspers. For a simple clasper
$C$, it is straightforward to show that there is a sequence of admissible
graph types
$s_0$ such that $s_0\geq (G)$, where $G$ is the graph type of $C$. Then we
have $(G)\leq s_1\leq s_2\leq \ldots \leq s_n\leq s_0$, where each sequence
is related to the next by a replacement of a single graph by a pair related
by IHX. 

 Proposition \ref{realize} 
implies that an elementary relation $s_i\leq s_{i+1}$ can be realized
geometrically modulo higher grope degree. More precisely, if
$s_i$ represents a sequence of clasper surgeries between a knot $K_0$ and
$K_1$, then there is a sequence of clasper surgeries of the form
$s_{i+1}$ which also go between $K_0$ and $K_1$. Therefore, we have that
$C$ can be implemented, modulo higher degree, by a sequence of claspers of the
form
$s_0$, which are admissible. This concludes the proof of part (a).

\noindent Part (b):
Since the clasper is admissible, every univalent vertex is part of a Y,
such that the two other vertices of the Y are trivalent, and there is a
path connecting them which doesn't hit the Y's interior.
 Surger the knot along the Y. This gives a connected clasper which has
two new capped leaves having two intersections with the knot each. Refining this
clasper using Theorem 20 of \cite{ct}, we will get some claspers of higher
grope degree and some simple claspers where some number of edges have been
cut. All of the latter type have increased
$c_1$, since we took a single simple cap of the original and converted it
to two simple caps in all the daughters. Cutting of edges will add even
more simple caps.
\end{proof}

\begin{lemma}
\begin{itemize}
\item[(a)]A clasper surgery of grope degree~$(2k-1)$ is realizable by a
sequence of simple tree clasper surgeries of grope degree~$(k-1)$.
\item[(b)]A clasper surgery which has $k+1$ caps is realizable by a sequence
of simple tree clasper surgeries of grope degree~$(k-1)$.
\end{itemize}
\end{lemma}
\begin{proof}The proofs of both of these facts use the main theorem of Habiro
\cite{h2}, which is that if two knots have the same degree~$(k-1)$
invariants, then they are related by a sequence of simple tree clasper
surgeries of Vassiliev=grope degree~$k$. The fact that grope degree~$(2k-1)$
clasper surgeries preserve type $k-1$ invariants is Theorem 3 of
\cite{newconant}. For part (b), if a rooted clasper $C$ has $k+1$ caps, it
must have $k$ non-root caps. Then there are $k$ groups of crossing changes
on $K_C$ which correspond to pushing the knot out of the caps. 
Since surgery on
a clasper which has a cap that does not intersect the knot is trivial, doing
any collection of these $k$ groups of crossings will yield the knot $K$. Hence
$K$ and $K_C$ are $k-1$-equivalent. 
\end{proof}

\begin{proof}[Proof of Theorem~\ref{imply}]
Define a complexity
function, ordered lexicographically as the triplet
$(c_1,g,a)$ where $c_1$ is the number of simple capped leaves, $g$ is the
grope degree,  and $a(G)$ is zero unless $G$ is admissible, in which case
it is one. If $g\geq 2k-1$ or $c_1\geq k+1$ then by the previous lemma the
surgery is realizable by simple grope degree
$k-1$ surgeries and we are done. 

We prove the following statement by contradiction: ``Every
clasper surgery of grope degree~$\geq k-1$ is realizable by a sequence of
simple clasper surgeries of grope degree~$(k-1)$."

Assume $C$ is a counterexample. Then, as noted above, $C$ lies inside the range
$c_1\leq k$ and
$k-1\leq g\leq 2k-2$. Hence it makes sense to take $C$ to
have maximal complexity $(c_1,g,a)$. There are two cases: either $g(C)=k-1$ or
$g(C)\geq k$. In the former case, $C$ cannot be admissible, since
it would then be simple and hence not a counterexample. Therefore, $a=0$, in
which case Lemma~\ref{lemmagraphmoves}(a) says that $C$ can be realized by
claspers with higher complexity. Since $C$ is a counterexample, one of these
daughters must also be a counterexample. But this contradicts that $C$ was of
maximal complexity.

So we are left with the case that $g(C)\geq k$. 
If $a=0$, then as above Lemma~\ref{lemmagraphmoves}(a) furnishes a
contradiction. Otherwise, $a=1$ and so
Lemma~\ref{lemmagraphmoves}(b) says that $C$ can be realized by claspers 
of degree~$\geq k-1$, and of higher complexity. As before, this contradicts the
maximality of $C$.
  \end{proof}


\begin{thebibliography}{999}
\bibitem[C]{c} J. Conant, \textsl{A knot bounding a grope of class
$n$ is $\lceil \frac{n}{2} \rceil$-trivial}, math.GT/9907158
\end{thebibliography}
\end{document}